\documentstyle{amsppt}
\baselineskip18pt
\magnification=\magstep1
\pagewidth{30pc}
\pageheight{45pc}

\hyphenation{co-deter-min-ant co-deter-min-ants pa-ra-met-rised
pre-print pro-pa-gat-ing pro-pa-gate
fel-low-ship Cox-et-er dis-trib-ut-ive}
\def\leaderfill{\leaders\hbox to 1em{\hss.\hss}\hfill}
\def\A{{\Cal A}}

\def\H{{\Cal H}}

\def\Pl{{\Cal P}}
\def\Sy{{\Cal S}}\

\def\afn{{\text {\bf a}}}
\def\tr{{\text {tr}}}

\def\idest{i.e.,\ }

\def\a{{\alpha}}

\def\g{{\gamma}}
\def\G{{\Gamma}}
\def\d{{\delta}}

\def\e{{\varepsilon}}

\def\k{{\kappa}}
\def\l{{\lambda}}

\def\s{{\sigma}}
\def\t{{\tau}}

\def\bc{{\bold c}}

\def\BB{{\bold B}}

\def\b0{\text{\bf 0}}

\def\ra{{\ \longrightarrow \ }}

\def\lan{{\langle}}
\def\ran{{\rangle}}

\def\real{{\Bbb R}}
\def\complex{{\Bbb C}}
\def\zed{{\Bbb Z}}
\def\kyu{{\Bbb Q}}
\def\enn{{\Bbb N}}

\def\Im{\text{\rm Im}}

\def\boxit#1{\vbox{\hrule\hbox{\vrule \kern3pt
\vbox{\kern3pt\hbox{#1}\kern3pt}\kern3pt\vrule}\hrule}}
\def\rabbit{\vbox{\hbox{\kern0pt
\vbox{\kern0pt{\hbox{---}}\kern3.5pt}}}}

\def\tableau#1{
        \hbox {
                \hskip -10pt plus0pt minus0pt
                \raise\baselineskip\hbox{
                \offinterlineskip
                \hbox{#1}}
                \hskip0.25em
        }
}

\def\tabCol#1{
\hbox{\vtop{\hrule
\halign{\strut\vrule\hskip0.5em##\hskip0.5em\hfill\vrule\cr\lower0pt
\hbox\bgroup$#1$\egroup \cr}
\hrule
} } \hskip -10.5pt plus0pt minus0pt}

\def\CR{
        $\egroup\cr
        \noalign{\hrule}
        \lower0pt\hbox\bgroup$
}



\def\blank#1#2{
\hbox to #1{\hfill \vbox to #2{\vfill}}
}


\def\strut{\vrule height10pt depth5pt width0pt}


\def\ignore#1{{}}

\topmatter
\title Constructing cell data for diagram algebras
\endtitle

\author R.M. Green and P.P. Martin \endauthor
\affil Department of Mathematics \\ University of Colorado \\
Campus Box 395 \\ Boulder, CO  80309-0395 \\ USA \\ {\it  E-mail:}
rmg\@euclid.colorado.edu \\
\newline
Centre for Mathematical Science \\ City University \\
Northampton Square \\ London EC1V 0HB \\ England \\ {\it E-mail:}
P.P.Martin\@city.ac.uk \\
\endaffil

\abstract 
We show how the treatment of cellularity in families of algebras arising
from diagram calculi, such as Jones' Temperley--Lieb wreaths,
variants on Brauer's centralizer algebras,
and the contour algebras of Cox et al (of which many algebras are
special cases), may be unified using the theory of tabular algebras.
This improves an earlier result of the first author (whose hypotheses
covered only the Brauer algebra from among these families).
%
\endabstract


\endtopmatter

\centerline{\bf To appear in the Journal of Pure and Applied Algebra}


\head Introduction \endhead

Cellular algebras were introduced by Graham and Lehrer \cite{{\bf 11}},
and are a class of finite dimensional associative algebras defined in
terms of a ``cell datum'' and three axioms.  The axioms allow one to
define a set of modules for the algebra known as ``cell modules'', and
one of the main strengths of the theory is that it is relatively 
straightforward to construct and to classify the irreducible modules 
for a cellular algebra in terms of quotients of the cell modules.

Tabular algebras were introduced by the first 
author in \cite{{\bf 14}} as a class of
associative $\zed[v, v^{-1}]$-algebras equipped with distinguished bases
(tabular bases) and satisfying certain axioms.  
In the most general setting, tabular algebras are defined via a somewhat
complicated ``table datum'' extending the cell datum construct.  However, 
there is a large natural subclass, the so-called ``tabular algebras with
trace'', which may be defined (up to isomorphism, in
a sense made precise in \cite{{\bf 16}}) simply by giving the distinguished basis.
In \cite{{\bf 17}}, the first author introduced ``cell modules'' and ``standard
modules'' for tabular algebras; each of these classes of modules is 
analogous in some sense to the cell modules of a cellular algebra.

The motivation for the theory of tabular algebras is twofold.  On one hand,
the theory can be viewed as a framework for studying the properties of 
``canonical bases'' for algebras.  The latter objects have been studied 
abstractly using constructions such as Du's IC bases \cite{{\bf 8}} and 
Stanley's $P$-kernels \cite{{\bf 32}}, and the archetypal examples are the
celebrated Kazhdan--Lusztig bases introduced in \cite{{\bf 21}}.  This is the
point of view taken in \cite{{\bf 17}}, where it is shown that the Kazhdan--Lusztig
bases of certain extended affine Hecke algebras of type $A$ are tabular
bases, and that the standard modules for the tabular algebra agree with
the geometrically defined standard modules appearing in the work of Lusztig 
\cite{{\bf 25}} and others.

The other main motivation for the theory of tabular algebras is as templates 
for cellular algebras.  Core to this is a theorem \cite{{\bf 14}, Theorem 2.1.1}
giving conditions under which one can describe a cellular structure for an 
algebra in
terms of the tabular structure.  As we will explain in this paper, there are
many cases in the literature where a cellular basis for an algebra has been
constructed in terms of another basis that turns out to be tabular; in these
situations, the tabular algebras may thus be regarded as more basic objects
than the corresponding cellular algebras.  This is a helpful point of view 
because tabular bases have some advantages over cellular 
bases: they are defined integrally (over $\zed[v, v^{-1}]$), and they are 
easy to construct in many cases because they occur naturally in the contexts of
Kazhdan--Lusztig type bases and algebras given by diagram calculi (many of 
which arise, for example, in computational statistical mechanics 
\cite{{\bf 26}, {\bf 30}, {\bf 34}}).

The main purpose of this paper is to generalize the core theorem 
\cite{{\bf 14}, Theorem 2.1.1} 
to cover a wider class of examples, with particular emphasis on examples 
that arise from algebras given by a calculus of diagrams.
We thus obtain shorter proofs of cellularity in several of the known
examples of cellular algebras (see for example corollaries 4.2.4, 5.2.3,
and 5.3.6): the constructions are similar, but our
approach has the advantage of relative generality.  One can also use our 
main result,
Theorem 2.2.1, to construct new examples of cellular algebras.  One way
to do this is using Theorem 4.1.3, which shows how to construct a kind
of wreath product of certain cellular algebras, which in turn yields new 
examples of
cellular algebras.  Another method we use involves the notion of a subdatum
(introduced in Definition 2.1.4), which is a convenient way to describe 
certain subalgebras of tabular algebras as tabular algebras in their own right
(Proposition 2.1.5).  It may be anticipated that these techniques will
provide further short proofs of cellularity in future applications.

Another approach to finding cellular bases for certain diagram algebras was
given by Enyang \cite{{\bf 9}}, who showed how to lift a cellular basis for the
Hecke algebra to a cellular basis for the Birman--Murakami--Wenzl algebra,
which is a $q$-analogue of the Brauer algebra.  It would be interesting
to know if Enyang's technique can be related to that of this paper.


\head Diagram Algebras and Cellular Algebras \endhead

\def\mbox#1{\text{\rm #1}}
\def\Hom{{\mbox{Hom}}}
\def\trungent{flush}

The formal definition of {\it diagram algebra} is beyond the scope of
this paper (see \cite{{\bf 29}} for the core paradigm),
but there are some simple components which it will be useful to bring
to mind, in order most simply to complete  discussion of 
the historical context of our work. 
 
Let $X$ be a poset. 
A {\it formal diagram category} (on $X$) is a category whose objects
are the elements of $X$, and morphisms $d \in \Hom(x,y)$ are called diagrams, 
with properties (some of) which are described below.
If $d \in \Hom(x,y)$ may be expressed as 
$d_1 d_2$ with $d_1 \in \Hom(x,x')$ and $d_2 \in \Hom(x',y)$
we say $d$ {\it factors} through $x'$.
A {\it propagating index} of $d$ is a lowest element of $X$ such that
$d$ factors through it. Such is not unique in general: we let 
$\# d$ denote the set of propagating indices of $d$. 
A diagram category has the {\it filtration property} 
$$
\# d_1 d_2 \leq \# d_1 , \; \# d_2
$$
i.e. $x \in \# d_1 d_2$ implies $x \leq y$ for $y \in \# d_i$.

If $d \in \Hom(x,x)$ has $x$ as a propagating index it is said to be
{\it \trungent}. In particular, $1_x$ is \trungent. 
The subset of \trungent\ diagrams is denoted $\Hom_t(x,x)$. 
In a diagram category the composite of \trungent\ diagrams is \trungent, 
so $\Gamma(x)=\Hom_t(x,x)$ is a kind of submonoid of $\Hom(x,x)$.

As an example, 
let $X$ be the set of natural numbers with the natural total order.
Then we may associate a category on $X$ in which the morphisms are 
Temperley--Lieb diagrams.
(These are defined formally later in the paper, but for now we
give a heuristic description.) 
A Temperley--Lieb diagram in $\Hom(m,n)$ 
is a set of $m+n$ vertices on the boundary 
of an interval of the plane,
together with a non-crossing partition into pairs
of these vertices.
{\it Non-crossing} is 
the property that the pairings
may be realised by connecting line
segments between the vertices, drawn in the 
interior of the interval without crossing each other. 
(Such a diagram is illustrated in Figure~1.)
Composition of morphisms may be computed by concatenation of 
diagrams so that the last $n$ vertices of $d_1$ meet 
the first $n$ of $d_2$ and become internal points. 
In this case there is only one propagating index, which is the
number of distinct lines which pass between the first $m$ and the last $n$
vertices in $d$. 
The only flush diagram in $\Hom(n,n)$ is $1_n$ itself. 
\vfill\eject
\topcaption{Figure~1} A diagram arising from a non-crossing partition,
with $m = n = 7$ and propagating index $1$ \endcaption
\centerline{
\hbox to 4.083in{
\vbox to 2.180in{\vfill
        \includegraphics{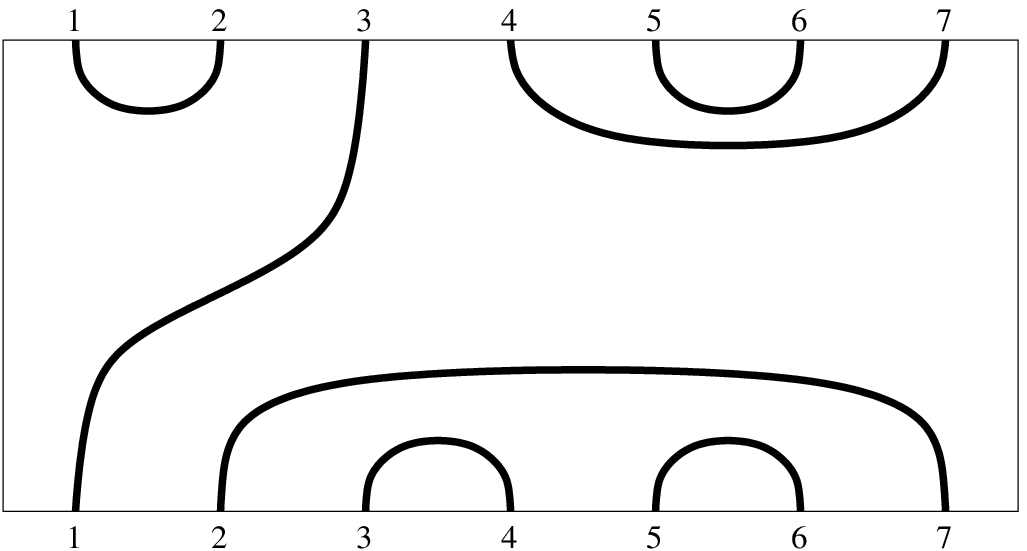}
}
\hfill}
}

Let us focus for a moment on the observation that $\Hom(x,y)$
has an action of $\Hom(x,x)$ on the left and $\Hom(y,y)$ on the right.
(So far this set has only a left semigroup action of $\Hom(x, x)$ rather
than a module structure, but it will be
convenient to adopt (bi)module terminology.)
This  $\Hom(x,y)$ can be partitioned into components with given
propagating index, and the parts with propagating indices less than 
a given index ($z$ say) form a sub-bimodule, 
by the filtration property. Denote the quotient  $\Hom^z(x,y)$.
Suppose that $x \leq y$. 
Then (in a diagram category) $\Hom^x(x,y)$ is isomorphic to a sum of copies of 
$\Hom_t(x,x)$ with respect to its left action. 
This paper concerns (from the diagram algebra perspective) the
inheritance of properties of  $\Hom_t(x,x)$--modules by 
the $x$-section of the above filtration
(again, see \cite{{\bf 29}} for concrete examples).

Diagram algebras and cellular algebras have a history of intertwined
development, and this inheritance aspect is no exception. 
The first several notable examples of cellular bases appeared {\it
before} the introduction of cellular algebras, in diagram algebra
contexts (see \cite{{\bf 18}, {\bf 27}, {\bf 30}}). 
In the physical contexts in which diagram algebras occur (e.g. as
transfer matrix algebras) the cellular axioms which we will describe
in the next section correspond heuristically, but closely, to ideas of
information propagation and of spatial or time-reversal symmetry. 
On the other hand, K\"onig and Xi \cite{{\bf 22}} have introduced
their ``inflationary'' construction for cellular algebras, which puts
these ideas in a nice abstract setting. 
A diagram algebra is an algebra with a ``natural'' basis of diagrams.
As noted above, each diagram may be cut into two parts: a top and a
bottom, say. The set of possible top parts is taken into the set of
bottom parts by an operation $i$ turning them upside down. 
This set is partitioned into subsets (``by) layers'' $a \sim b$ if $a$,
$i(b)$ can be the parts of a cut diagram.
However the combination of $a$, $i(b)$ is not in general unique, this
combination being controlled by an intermediate algebra. 
Thus layers of the algebra take the form $V \otimes G \otimes V$
where $G$ is the intermediate. 

K\"onig and Xi noted that subject to some technical conditions
consistent with the above, a rather general free $R$--module $V$
(irrespective of diagrams) may be used with a cellular algebra $G$ to
produce another cellular algebra with layer $V \otimes G \otimes V$. 
The argument hinges on an equivalent definition of cellular algebra
that does not use bases, and the construction does not provide any
``inflated'' basis. 

In practical matters of representation theory of concrete algebras,
however, explicit cellular (indeed any) bases are invaluable. 
Here, using the first author's relatively robust ``tabular bases'' and 
algebras $G$
that are hypergroups we are able to develop a version of inflation
with bases. As we will see, these bases can be chosen to {\it be} the
natural bases in diagram algebra examples. 
Thus our theorem offers a formalism well tuned to analysis of such
concrete algebras, a useful counterpoint to K\"onig and Xi's
elegantly abstract construction. 

\head 1. Cellular algebras and hypergroups \endhead

\subhead 1.1 Cellular algebras \endsubhead

Cellular algebras were originally defined by Graham and Lehrer
\cite{{\bf 11}}.

\definition{Definition 1.1.1}
Let $R$ be a commutative
ring with identity.  A {\it cellular algebra} over $R$ is an associative unital
algebra, $A$, together with a cell datum $(\Lambda, M, C, *)$ where:

\item {(C1)}
{$\Lambda$ is a finite poset.  For each $\l \in \Lambda$, $M(\l)$ is a 
finite set such that $$
C : \coprod_{\l \in \Lambda} \left( M(\l) \times M(\l) \right) \rightarrow A
$$ is injective with image an $R$-basis of $A$.}
\item {(C2)}
{If $\l \in \Lambda$ and $S, T \in M(\l)$, we write $C(S, T) = C_{S, T}^{\l}
\in A$.  Then $*$ is an $R$-linear involutory anti-automorphism 
of $A$ such that
$(C_{S, T}^{\l})^* = C_{T, S}^{\l}$.}
\item {(C3)}
{If $\l \in \Lambda$ and $S, T \in M(\l)$ then for all $a \in A$ we have $$
a . C_{S, T}^{\l} \equiv \sum_{S' \in M(\l)} r_a (S', S) C_{S', T}^{\l}
\mod A(< \l),
$$ where  $r_a (S', S) \in R$ is independent of $T$ and $A(< \l)$ is the
$R$-submodule of $A$ generated by the set $$
\{ C_{S'', T''}^{\mu} : \mu < \l, S'' \in M(\mu), T'' \in M(\mu) \}
.$$}
\enddefinition

\remark{Remark 1.1.2}
We have assumed $\Lambda$ to be finite to avoid complications (see
\cite{{\bf 12}, \S1.2}).
\endremark

We now recall from the literature some of the main examples of cellular of 
algebras that are particularly relevant for our purposes in this paper.

\example{Example 1.1.3}
Let ${\Cal S}_n$ be the symmetric group on $n$ letters.  Then the group
algebra $\zed {\Cal S}_n$ is cellular over $\zed$.  In this case, the
poset $\Lambda$ is the set of partitions of $n$, ordered by dominance (meaning
that if $\l \trianglerighteq \mu$ then $\l \leq \mu$).  The set $M(\l)$
is the set of standard tableaux of shape $\l$, namely the ways of writing
the numbers $1, \ldots, n$ once each into a Young diagram of shape $\lambda$
such that the entries increase along rows and down columns.  The element
$C_{S, T}^\lambda$ is the Kazhdan--Lusztig basis element $C_w$ such that
$w \in {\Cal S}_n$ corresponds via the Robinson--Schensted correspondence to
the ordered pair of standard tableaux $(S, T)$.  The map $*$ sends $C_w$
to $C_{w^{-1}}$.

The Hecke algebra $\H({\Cal S}_n)$ was shown to be cellular by Graham and
Lehrer in \cite{{\bf 11}, Example 1.2}, and the underlying idea was already 
implicit in \cite{{\bf 21}}.  The example of the symmetric group above is 
obtained simply 
by specializing $q$ to $1$, as was observed by Graham and Lehrer in their
treatment of the Brauer algebra \cite{{\bf 11}, \S4}
For details on the relationship between
the Robinson--Schensted correspondence and Kazhdan--Lusztig theory, 
the reader is referred to Ariki's paper \cite{{\bf 3}}.
\endexample

\example{Example 1.1.4}
A simple example of a cellular algebra that is important for our purposes
is Graham and Lehrer's so-called ``banal example'' \cite{{\bf 11}, Example 1.3}.  
Let $R$ be a commutative ring with identity, let $\l_1, \l_2, \ldots, \l_k$ be (not necessarily
distinct) elements of $R$, and let $P(x) \in R[x]$ be the polynomial $
\prod_{i = 1}^k (x - \l_i)
.$  Then the rank $k$ algebra $A = R[x]/\lan P(x) \ran$ is cellular
over $R$.  A cell datum is as follows:  $\Lambda$ is the poset 
$\{1, 2, \ldots, k\}$, ordered
in the natural way, $M(\lambda)$ is a one-element set for each $\lambda$ and
$C_{S, S}^j$ is the image of the polynomial $
\prod_{i = j + 1}^k (x - \l_i)
.$  The map $*$ is the identity map, which is an anti-automorphism because
$A$ is commutative.
\endexample

Let $A$ and $A'$ be cellular algebras over $R$ with cell data
$(\Lambda_1, M_1, C_1, *_1)$ and $(\Lambda_2, M_2, C_2, *_2)$ respectively.
We will show in the next two examples how the direct sum and direct product 
of two cellular algebras are again cellular in a natural way.
(We omit the proofs of these results because they are both well-known and easy:
see the remarks at the end of \cite{{\bf 22}, \S6}.)
Later (Theorem 4.1.3) we will look at a less trivial way to form new cellular
algebras using a kind of wreath product.

\example{Example 1.1.5}
Let $\Lambda_3$ be the disjoint union of $\Lambda_1$ and $\Lambda_2$.
We partially order order $\Lambda_3$ by stipulating that $\l \leq \l'$ if,
for some $i \in \{1, 2\}$, we have $\l, \l' \in \Lambda_i$ and $\l \leq_i \l'$,
where $\leq_i$ is the partial order on $\Lambda_i$.
For $\l \in \Lambda_3$, we define $M_3(\l)$ to be
$M_1(\l)$ if $\l \in \Lambda_1$ and $M_2(\l)$ if $\l \in \Lambda_2$.
We define $C_3$ (respectively, $*_3$) in an analogous way as a natural 
extension of $C_1$ and $C_2$ (respectively, $*_1$ and $*_2$).
The $R$-algebra $A \oplus A'$ is then cellular with cell datum $(\Lambda_3, 
M_3, C_3, *_3)$.
\endexample

\example{Example 1.1.6}
Let $\Lambda_4$ be the Cartesian product $\Lambda_1 \times \Lambda_2$,
partially ordered by  stipulating that $(\l_1, \l_2) \leq_4 (\l'_1, \l'_2)$ 
if and only if $\l_1 \leq_1 \l'_1$ and $\l_2 \leq_2 \l'_2$.
For $\l = (\l_1, \l_2) \in \Lambda_4$, we define $M_4(\l)$ to be
$M_1(\l_1) \times M_2(\l_2)$.  For $(S_1, S_2), (T_1, T_2) \in M(\l)$, we 
define $C_4((S_1, S_2), (T_1, T_2))$ to be $C_1(S_1, T_1) \otimes_R
C_2(S_2, T_2)$.  The map $*_4$ is the $R$-linear map sending
$C_4((S_1, S_2), (T_1, T_2))$ to $C_4((T_1, T_2), (S_1, S_2))$.
The $R$-algebra $A \otimes_R A'$ is then cellular with cell datum $(\Lambda_4, 
M_4, C_4, *_4)$.
\endexample

\subhead 1.2 Hypergroups \endsubhead

A key ingredient in the definition of tabular algebras is the notion of a 
hypergroup.  There are many variants of this idea in the literature, for
example the table algebras of Arad and Blau \cite{{\bf 1}}, the generalized
table algebras of Arad, Fisman and Muzychuk \cite{{\bf 2}}, the association
schemes of algebraic combinatorics \cite{{\bf 4}} and the discrete
hypergroups as described by Sunder \cite{{\bf 33}, Definition IV.1}.  The
hypergroups we define here are the ``normalized table algebras'' of
\cite{{\bf 14}, Definition 1.1.2}, but we use the name ``hypergroups'' here for 
simplicity and to reflect the fact that most of the important examples
of hypergroups considered in this paper are in fact groups.

\definition{Definition 1.2.1}
A hypergroup is a pair $(A, B)$, where $A$ is an associative 
unital $R$-algebra for some commutative ring $R$ with $1$ and 
containing $\zed$,
and $B = \{b_i : i \in I\}$ is a distinguished basis for
$A$ such that $1 \in B$, satisfying the following three axioms:

\item{(H1)}{The structure constants of $A$ with respect to the basis
$B$ are nonnegative integers.}
\item{(H2)}{There is an algebra anti-automorphism $\bar{\ }$ of $A$ whose
square is the identity and that has the property 
that $b_i \in B \Rightarrow \overline{b_i} \in
B$.  (We define $\overline{i}$ by the condition $\overline{b_i} =
b_{\bar{i}}$.)}
\item{(H3)}{Let $\k(b_i, a)$ be the coefficient of $b_i$ in $a \in A$.
Then we have $\k(b_m, b_i b_j) = \k(b_i, b_m \overline{b_j})$
for all $i, j, m \in B$.}
\enddefinition

\remark{Remark 1.2.2}
Note that, setting $b_m = 1$ in axiom (H3), we see that $b_{\bar{i}}$ can
be characterized by the property that it is the unique basis element $b_j$ for
which $1$ appears with nonzero coefficient in $b_i b_j$ (or in $b_j b_i$).
This implies that the anti-automorphism $\bar{\ }$ is completely determined
by the structure constants.
\endremark

\example{Example 1.2.3}
Perhaps the most obvious example of a hypergroup is the case where $B$
is a group $G$ and $A$ is the group algebra $RG$, where $R$ is a commutative
ring with 1 that contains $\zed$.  In this case, the anti-automorphism 
$\bar{\ }$ is the $R$-linear extension of inversion in $G$.
\endexample

The following result is an easy consequence of axiom (H3) (see also 
\cite{{\bf 15}, Proposition 1.1.4}).

\proclaim{Proposition 1.2.4}
Let $(A, B)$ be a hypergroup.  The linear function $t$ 
sending $a \in A$ to $\k(1, a)$ satisfies $t(xy) = t(yx)$ for all $a \in A$.
\qed\endproclaim

The following result shows how a tensor product of two hypergroups is another
hypergroup.  (In the case of groups, this construction corresponds to the
direct product.)

\proclaim{Proposition 1.2.5}
Let $(A_1, B_1)$ and $(A_2, B_2)$ be hypergroups over $R$.
Then $$(A_1 \otimes_R A_2, B_1 \otimes_R B_2)$$ is a hypergroup over $R$,
where the multiplication on $A_1 \otimes A_2$ is given by the
Kronecker product and the anti-automorphism $\bar{\
}$ of $A_1 \otimes_R A_2$ is defined to send $b_1 \otimes_R b_2$ to
$\overline{b_1} \otimes_R \overline{b_2}$.
\endproclaim

\demo{Proof}
See, for example, \cite{{\bf 15}, Proposition 1.1.5}.
\qed\enddemo

\subhead 1.3 Based rings \endsubhead

We now recall Lusztig's notion of a based ring; see \cite{{\bf 24}} or 
\cite{{\bf 37}, \S1.5}.

\definition{Definition 1.3.1}
A {\it based ring} is a pair $(A, B)$, where $A$ is a unital $\zed$-algebra
with free $\zed$-basis $B$ and nonnegative structure constants.  A
homomorphism $\phi : (A, B) \ra (A', B')$ of based rings is a
homomorphism of abstract $\zed$-algebras $\phi : A \ra A'$ such that $\phi(b)
\in B' \cup \{0\}$ for all $b \in B$.  Isomorphisms,
automorphisms, anti-automorphisms, etc. of based rings are defined 
analogously.  
\enddefinition

\remark{Remark 1.3.2}
Clearly hypergroups over $\zed$ give examples of based rings, and the map 
$\bar{\ }$ of axiom (H2) is an anti-automorphism of based rings.
\endremark

\proclaim{Proposition 1.3.3}
Let $(A, \BB)$ and $(A', \BB')$ be hypergroups over $\zed$ and let 
$$f : (A, \BB) \ra (A', \BB')$$ be a unital homomorphism of based rings.  Then,
for each $a \in A$, we have $f(\bar{a}) = \overline{f(a)}$, where the 
anti-automorphisms $\bar{\ }$ are those associated by Definition 1.2.1
to each hypergroup.
\endproclaim

\demo{Proof}
Because $f$ is $\zed$-linear, we may immediately reduce our consideration 
to the case where $a \in \BB$.

Let $b_i \in \BB$, and consider the equation $$
b_i b_{\bar i} = 1 + \sum_{1 \ne b_k \in \BB} c_k b_k
,$$ which holds by Remark 1.2.2.  Applying
$f$ to the equation and using the hypothesis that $f$ sends
$1 \in \BB$ to $1 \in \BB'$, we obtain $$
f(b_i) f(b_{\bar i}) = 1 + \sum_{1 \ne b_k \in \BB} c_k f(b_k)
.$$   Although it may be the case that $f(b_k) = 1$ for some $b_k \ne 1$,
axiom (H1) shows that the coefficient of $1$ on the right hand side is
nonzero.  Remark 1.2.2 now shows that $f(b_{\bar i}) = f(\overline{b_i}) = 
\overline{f(b_i)}$, where the first equality is by definition of 
$b_{\bar i}$.  
\qed\enddemo

Proposition 1.3.3 ensures that the following definition makes sense.

\definition{Definition 1.3.4}
Let $(A, \BB)$ be a hypergroup over $\zed$ and let $$f : (A, \BB) \ra
(A, \BB)$$ be an automorphism of based rings.  Then the anti-automorphism
$$\overline{f} : (A, \BB) \ra (A, \BB)$$ is defined to be the composition
of $f$ with the anti-automorphism $\bar{\ }$ of axiom (H2).
\enddefinition

\proclaim{Proposition 1.3.5}
Let $(A, \BB)$ be a hypergroup over $\zed$ and let $$f : (A, \BB) \ra
(A, \BB)$$ be an anti-automorphism of based rings.  Then $f$ is of the
form $\overline{g}$ for a unique automorphism $g : (A, \BB) \ra
(A, \BB)$ of based rings.
\endproclaim

\demo{Proof}
Composing $f$ with the hypergroup anti-automorphism $\bar{\ }$, we obtain
an automorphism $g$ with the required properties.  The uniqueness follows
from the invertibility of $\bar{\ }$.
\qed\enddemo

\head 2. Tabular algebras \endhead

We now recall the definition of tabular algebras from \cite{{\bf 14}}.

\subhead 2.1 Definition \endsubhead

\definition{Definition 2.1.1}
Let $\A = \zed[v, v^{-1}]$.  A {\it tabular algebra} is an
$\A$-algebra $A$, together with a table datum 
$(\Lambda, \Gamma, B, M, C, *)$ where:

\item{(A1)}
{$\Lambda$ is a finite poset.  For each $\l \in \Lambda$, 
$(\Gamma(\l), B(\l))$ is a hypergroup over $\zed$ and
$M(\l)$ is a finite set.  The map $$
C : \coprod_{\l \in \Lambda} \left( M(\l) \times B(\l) \times M(\l)
\right) \rightarrow A
$$ is injective with image an $\A$-basis of $A$.  We assume
that $\Im(C)$ contains a set of mutually orthogonal idempotents 
$\{1_\e : \e \in {\Cal E}\}$ such that 
$A = \sum_{\e, \e' \in {\Cal E}} (1_\e A 1_{\e'})$ and such that for each
$X \in \Im(C)$, we have $X = 1_\e X 1_{\e'}$ for some $\e, \e' \in
{\Cal E}$.
A basis arising in
this way is called a {\it tabular basis}.  
}
\item{(A2)}
{If $\l \in \Lambda$, $S, T \in M(\l)$ and $b \in B(\l)$, we write
$C(S, b, T) = C_{S, T}^{b} \in A$.  
Then $*$ is an $\A$-linear involutory anti-automorphism 
of $A$ such that
$(C_{S, T}^{b})^* = C_{T, S}^{\overline{b}}$, where $\bar{\ }$ is the
hypergroup anti-automorphism of $(\Gamma(\l), B(\l))$.
If $g \in \complex(v) \otimes_\zed \Gamma(\l)$ is such that 
$g = \sum_{b_i \in B(\l)} c_i b_i$ for some scalars $c_i$ 
(possibly involving $v$), we write
$C_{S, T}^g \in \complex(v)\otimes_\A A$ 
as shorthand for $\sum_{b_i \in B(\l)} c_i C_{S, T}^{b_i}$.  We write
$\bc_\l$ for the image under $C$ of $M(\l) \times B(\l) \times M(\l)$.}
\item{(A3)}
{If $\l \in \Lambda$, $g \in \Gamma(\l)$ and $S, T \in M(\l)$ then for all 
$a \in A$ we have $$
a . C_{S, T}^{g} \equiv \sum_{S' \in M(\l)} C_{S', T}^{r_a(S', S) g}
\mod A(< \l),
$$ where  $r_a (S', S) \in \Gamma(\l)[v, v^{-1}] = \A \otimes_\zed
\Gamma(\l)$ is independent of $T$ and of $g$ and $A(< \l)$ is the
$\A$-submodule of $A$ generated by the set $\bigcup_{\mu < \l} \bc_\mu$.}
\enddefinition

In all the examples we consider in this paper, the tabular basis will contain
the identity element of the algebra.  This means that the set ${\Cal E}$ 
contains only the identity element of $A$.

The paper \cite{{\bf 14}} also defines a more restrictive class of tabular
algebras called ``tabular algebras with trace''.  Since we are mainly
concerned with representation theory and not Kazhdan--Lusztig theory in this
paper, tabular algebras with trace will not be our primary objects of study.
However, we recall the definition here for later reference.
To do this, we need to recall the notion of $\afn$-function, due to 
Lusztig.

\definition{Definition 2.1.2}
Let $g_{X, Y, Z} \in \A$ be one of the structure constants for the tabular
basis $\Im(C)$ of $A$, namely $$
X Y = \sum_Z g_{X, Y, Z} Z
,$$ where $X, Y, Z \in \Im(C)$.   Define, for $Z \in \Im(C)$, $$
\afn(Z) = 
\max_{X, Y \in \Im(C)} \deg(g_{X, Y, Z})
,$$ where the degree of a Laurent polynomial is taken to be 
the highest power of $v$ occurring with nonzero coefficient.  We
define $\g_{X, Y, Z} \in \zed$ to be the coefficient of $v^{\afn(Z)}$ in $g_{X,
Y, Z}$; this will be zero if the bound is not achieved.
\enddefinition

\definition{Definition 2.1.3}
A {\it tabular algebra with trace} is a tabular algebra in the sense
of Definition 2.1.1 that satisfies the conditions (A4) and (A5) below.

\item{(A4)}{Let $K = C_{S, T}^b$, $K' = C_{U, V}^{b'}$ and 
$K'' = C_{X, Y}^{b''}$ lie in $\Im(C)$.  Then the
maximum bound for $\deg(g_{K, K', K''})$
in Definition 2.1.2 is achieved if and only if $X = S$, $T = U$, $Y =
V$ and $b''$ occurs with nonzero coefficient in $bb'$.
If these conditions all hold and
furthermore $b = b' = b'' = 1$, we require $\g_{K, K', K''} = 1$.}
\item{(A5)}{There exists an $\A$-linear function $\t : A \ra \A$
(the {\it tabular trace}), such that $\t(x) = \t(x^*)$ for all $x \in
A$ and $\t(xy) = \t(yx)$ for all $x, y \in A$, that has the 
property that for every
$\l \in \Lambda$, $S, T \in M(\l)$, $b \in B(\l)$ and $X = C_{S,
T}^b$, we have $$
\t(v^{\afn(X)} X) = 
\cases 1 \mod v^{-1} \A^- & \text{ if } S = T \text{ and } b = 1,\cr
0 \mod v^{-1} \A^- & \text{ otherwise.} \cr
\endcases
$$  Here, $\A^- := \zed[v^{-1}]$.}
\enddefinition

The following notion is convenient for describing certain tabular algebras 
as subalgebras of other tabular algebras appearing in this paper.

\definition{Definition 2.1.4}
Let $A$ be a tabular algebra with table datum $(\Lambda, \Gamma, B, M, C, *)$.
A {\it subdatum} of such a table datum is a tuple
$(\Lambda', \Gamma', B', M', C', *')$ 
such that:

\item{(S1)}{$(\Lambda', \leq')$ is a subposet of $(\Lambda, \leq)$;}
\item{(S2)}{for each $\l' \in \Lambda'$, $M'(\l')$ is a subset of $M(\l')$
and there is a unital monomorphism of based
rings $(\Gamma'(\l'), B'(\l')) \ra (\Gamma(\l'), B(\l'))$ identifying
$(\Gamma'(\l'), B'(\l'))$ with a subhypergroup of $(\Gamma(\l'), B(\l'))$;}
\item{(S3)}{under the above identifications, the maps $C'$ and $*'$ are 
the restrictions of $C$ and $*$, respectively, and $\Im(C') = A'$ is an
$\A$-subalgebra of $A$.}
\enddefinition

The tuple defined above turns out to be a table datum for $A'$, as we now
show.

\proclaim{Proposition 2.1.5}
Let $A$ be a tabular algebra with table datum $(\Lambda, \Gamma, B, M, C, *)$,
and let $(\Lambda', \Gamma', B', M', C', *')$ be a subdatum for an 
$\A$-subalgebra $A'$ of $A$.  If the algebra $A'$ contains all the
idempotents $\{1_\e : \e \in {\Cal E}\}$ of axiom (A1) then  the given
subdatum is a table datum for $A'$.  If, 
furthermore, $A$ is a tabular algebra with trace, then so is $A'$.
\endproclaim

\demo{Proof}
We check the tabular axioms applied to $A'$.
Axiom (A1) follows from the definitions and the hypothesis about the
idempotents.  Axiom (A2) follows from the definitions and 
Proposition 1.3.3.  Axiom (A3) is immediate.

Suppose now that $A$ is a tabular algebra with trace.

Let $X = C_{S, T}^{\prime b} \in \Im(C')$, where $S, T \in M'(\l')$ and $b \in 
B'(\l')$ for some $\l' \in \Lambda'$.  
We first show that the $\afn$-functions $\afn_{A'}$ and
$\afn_A$ arising from the algebras $A'$ and $A$ take the same value on $X$.  
It is clear from the definition of the $\afn$-function and the fact that 
$\Im(C') \subseteq \Im(C)$ that $\afn_{A'}(X) \leq \afn_A(X)$.  For the
reverse inequality, we recall from \cite{{\bf 14}, Lemma 2.2.3} that $C_{S, T}^b$
occurs in the product $C_{S, S}^1 C_{S, T}^b$ with coefficient of degree
$\afn_A(C_{S, T}^b)$.  Now $C_{S, T}^b \in \Im(C')$ by definition of $X$,
and $C_{S, S}^1 \in \Im(C')$ because $S \in M'(\l')$ and $1 \in B'(\l')$
by the unital requirement of axiom (S2) in Definition 2.1.4.  Thus
$X$ occurs in the product $C_{S, S}^{\prime 1} C_{S, T}^{\prime b}$ 
with coefficient of degree
$\afn_A(X)$, and this implies that $\afn_{A'}(X) \geq \afn_A(X)$, as required.

Axiom (A4) follows from the aforementioned compatibility of 
$\afn$-functions and the definitions, and 
axiom (A5) follows by restricting the trace $\t$ of $A$ to $A'$.
\qed\enddemo


\subhead 2.2 The main result \endsubhead


We are now ready to state our main result.
Most of the rest of the paper will be devoted to studying examples of
Theorem 2.2.1.

\proclaim{Theorem 2.2.1}
Let $A$ be a tabular algebra of finite rank with table datum \newline
$(\Lambda, \Gamma, B, M, C, *)$; that is $|B(\l)| < \infty$ for each
$\l \in \Lambda$.

Let $R$ be a commutative ring with identity.
Suppose that $\a$ is an $R$-algebra automorphism of $A$ satisfying
$\a(\Im(C)) = \Im(C)$ and  with the property that, for each $\l \in \Lambda$, 
there exists a permutation $\s_\l$ of $M(\l)$ and an automorphism $f_\l$
of the based ring $(\Gamma(\l), B(\l))$ such that $$
\a(C_{S, T}^b) = C_{\s_\l(S), \s_\l(T)}^{f_\l(b)}
$$ for each $S, T \in M(\l)$ and $b \in B(\l)$.

Suppose that, for some $R \geq \zed$ and for each
$\l \in \Lambda$, the algebra $R \otimes_\zed \Gamma(\l)$ is cellular 
over $R$ with cell datum $(\Lambda_\l, M_\l, C_\l, \overline{f_\l})$, where 
$\overline{f_\l}$ is as in Definition 1.3.4.

Then $R \otimes_\zed A$ is cellular over $R \otimes_\zed \A$ with cell datum 
$(\Lambda', M', C', *')$, 
where $\Lambda' := \{(\l, \l') : \l \in \Lambda, \l' \in \Lambda_\l\}$ (ordered
lexicographically), $M'((\l, \l')) := M(\l) \times M_\l(\l')$, $C'((S,
s), (T, t))$ (where $(S, s), (T, t) \in M(\l) \times M_\l(\l')$) is
equal to $C_{S, \s_\l(T)}^{C_\l(s, t)}$ and $*' = * \circ \a = \a \circ *$,
so that $*' : C_{S, T}^b \mapsto C_{\s_\l(T), \s_\l(S)}^{\overline{f_\l}(b)}$.
\endproclaim

\demo{Proof}
Axiom (C1) for $R \otimes_\zed A$ follows from axiom (A1) applied to $A$ 
and axiom (C1) applied to each hypergroup $(\Gamma(\l), B(\l))$.

We have $* \circ \a = \a \circ *$ by Proposition 1.3.3.  It then follows from
axiom (A2) that $*' = * \circ \a = \a \circ *$ is an anti-automorphism, and
axiom (C2) follows because each hypergroup $R \otimes_\zed \Gamma(\l)$ 
is cellular with respect to the hypergroup anti-automorphism $\overline{f_\l}$.

To prove axiom (C3), let $\l \in \Lambda$ and 
let $C_\l(s, t)$ be a basis element of $\Gamma(\l)$ with $s, t \in
M_\l(\l')$.  Then by axiom (A3) we have, for any $a \in A$, $$
a . C_{S, \s_\l(T)}^{C_\l(s, t)} \equiv \sum_{S' \in M(\l)} C_{S',
\s_\l(T)}^{r_a(S', S) C_\l(s, t)} \mod A(< \l).
$$  Since $R \otimes_\zed \Gamma(\l)$ is cellular over $R$ with cell 
basis given by $C_\l$, it
follows by axiom (C3) applied to $R \otimes_\zed \Gamma(\l)$ that $$
r_a(S', S) C_\l(s, t) \equiv \sum_{s' \in M_\l(\l')} r'(S', S, s', s) 
C_\l(s', t)
\mod R \otimes \Gamma_\l (< \l'),
$$ where the $r'(S', S, s', s)$ are elements of $R \otimes_\zed \A$ that
are independent of $t$ (and, by axiom (A3), independent of $\s_\l(T)$).
Axiom (C3) follows by tensoring over $R$.
\qed\enddemo

\remark{Remark 2.2.2}
In the special case where the automorphism $\a$ is the identity map, Theorem
2.2.1 reduces to \cite{{\bf 14}, Theorem 2.1.1}.
\endremark

\remark{Remark 2.2.3}
The theorem can be proved using a weaker order on $\Lambda'$, 
namely the order such that $(\l_1, \l'_1) \leq (\l_2, \l'_2)$ if and only
if $\l_1 \leq \l_2$ and $\l'_1 \leq \l'_2$.  The proof is the same.
\endremark

The next result shows that if $\a$ is not the identity map, then it must be
an involution.

\proclaim{Lemma 2.2.4}
Let $\a$ be an automorphism satisfying the hypotheses of Theorem 2.2.1.  Then
$\a^2$ is the identity map.
\endproclaim

\demo{Proof}
This is immediate from the assumptions that $\a \circ *$ has order $2$,
$*$ has order $2$, and $\a$ commutes with $*$.
\qed\enddemo

It will also turn out that the map $\a$ in Theorem 2.2.1 need not be unique;
see Remark 5.3.7 below.

\remark{Remark 2.2.5}
The above results suggest a place to look for cellular involutions of a
tabular algebra $A$ in the case where the tabular involution does not work,
namely to compose the tabular involution with basis-preserving algebra 
automorphisms of order 2.
\endremark

\vfill\eject
\head 3. Diagram algebra preliminaries \endhead

In the rest of the paper 
we study examples of cellular algebras   
with bases consisting of diagrams of various kinds.
All of our examples can be related to each other.
Although the ordinary
Temperley--Lieb algebra (see \S3.2) is arguably the hub of these connections, 
it is convenient to start by recalling Brauer's centralizer algebra. 
Some useful references on this algebra are 
Brauer's original paper \cite{{\bf 6}}, as well as \cite{{\bf 11}, \S4} and 
\cite{{\bf 35}}.  


\subhead 3.1 Brauer diagrams \endsubhead

Combinatorially, the Brauer algebra $B_n$ has defining basis
consisting simply of the set of partitions of $2n$ objects into pairs.
It is natural, however, to provide a graphical realisation. 
We start by 
recalling Jones' formalism of $k$-boxes \cite{{\bf 20}}.  
For further details 
and references, the reader is referred to \cite{{\bf 15}, \S2}.

\definition{Definition 3.1.1}
Let $k$ be a nonnegative integer.  The {\it standard $k$-box}, ${\Cal
B}_k$, is the set $\{(x, y) \in \real^2 : 0 \leq x \leq k + 1, \ 0
\leq y \leq 1\}$, together with the $2k$ marked points $$\eqalign{
&1 = (1, 1), \ 
2 = (2, 1), \ 
3 = (3, 1), \ 
\ldots, \ k = (k, 1), \cr
&k + 1 = (k, 0), \
k + 2 = (k-1, 0), \
\ldots, \ 
2k = (1, 0).\cr
}$$
(This is called the {\it Temperley--Lieb numbering}. 
The {\it Brauer numbering} renumbers the points $k+i$ of the standard
$k$-box (for $1\leq i \leq k$) as $k+1-i$. See Figure~2.)  
\enddefinition

\definition{Definition 3.1.2}
Let $X$ and $Y$ be embeddings of some topological spaces (such as lines) 
into the standard $k$-box. 
Multiplication of such embeddings to obtain a new embedding in the
standard $k$-box shall, where appropriate, be defined via the following
procedure on $k$-boxes.  
The product $XY$ is the 
embedding
obtained by placing $X$ on top of $Y$ 
(that is, $X$ is first shifted in the plane by $(0,1)$ relative to $Y$, 
so that marked point $(i,0)$ in $X$ coincides with $(i,1)$ in $Y$),
rescaling vertically by a scalar factor of $1/2$ and applying the
appropriate translation to recover a standard $k$-box.
\enddefinition



\topcaption{Figure~2} A Brauer algebra basis element for $n = 6$ \endcaption
\centerline{
\hbox to 3.500in{
\vbox to 2.180in{\vfill
        \includegraphics{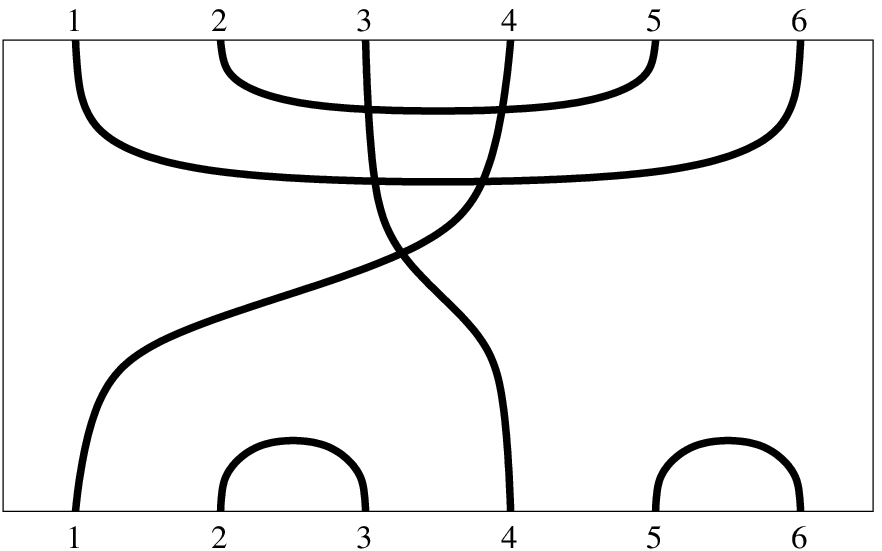}
}
\hfill}
}


\definition{Definition 3.1.3}
Let $k$ be a nonnegative integer. Consider the set of 
smooth embeddings of a single curve 
(which we usually call an ``edge'') in the standard $k$-box,
such that the curve is either closed (isotopic to a circle) 
or its endpoints coincide with two marked points of the box,
with the curve meeting the boundary of the box only at such
points, and there transversely.

By a smooth diffeomorphism of this curve we mean a smooth
diffeomorphism of the copy of $\real^2$ in which it is embedded,
that fixes the boundary, and in particular the marked points, of the
$k$-box, and takes the curve to another such smooth embedding. 
(Thus, 
the orbit of smooth diffeomorphisms
of one embedding contains all embeddings with the same endpoints.)

A concrete Brauer diagram is a set of such embedded curves with the property
that every marked point coincides with an endpoint of precisely one
curve. 
(In examples we can represent this set by drawing all the curves on
one copy of the $k$-box. Examples can always be chosen in which no 
ambiguity arises thereby; see Figure~2.)

Two such concrete diagrams are said to be equivalent if one may be
taken into the other 
by applying smooth diffeomorphisms to the individual curve
embeddings within it.

There is an obvious map from the set of concrete diagrams to the set
of pair partitions of the $2k$ marked points. 
It will be evident that the image under this map is an invariant of
concrete diagram equivalence.  

The set $B_k(\emptyset)$ is the set of equivalence classes of
concrete diagrams. Such a class (or any representative) is called
a Brauer diagram.  


Let $D_1, D_2$ be concrete diagrams. 
Since the $k$-box multiplication defined above
internalises marked points in coincident pairs, 
corresponding curve endpoints in $D_1D_2$ may also be internalised seamlessly. 
Each chain of curves concatenated in this way may thus be put in
natural correspondence with a single curve. 
Thus the multiplication 
gives rise to a closed associative binary operation on the
set of concrete diagrams.
It will be evident that this passes to a well defined
multiplication on $B_k(\emptyset)$. 
Let $R$ be a commutative ring with $1$.  
The elements of $B_n(\emptyset)$ form the basis elements of an $R$-algebra
$\Pl_n^B(\emptyset)$ with this multiplication. 

\enddefinition

A curve in a diagram   
that is not a closed loop
is called {\it propagating} if its endpoints have different
$y$-values, and {\it non-propagating} otherwise.  (Some authors use the
terms ``through strings'' and ``arcs'' respectively for curves of these
types.)

\definition{Definition 3.1.4}
The Brauer algebra $B_n = B_n(\d)$
is the free $R[\d]$-module with basis given by the elements
of $B_n(\emptyset)$ with no closed loops.  The multiplication is
inherited from the multiplication on $\Pl^B_n(\emptyset)$ except that
one multiplies by a factor of $\d = v + v^{-1}$ for each resulting closed 
loop and then discards the loop.

\enddefinition

In \S5 we shall return to consider the tabularity of $B_n$ and various
related algebras.  The assumption $\d = v + v^{-1}$ in Definition 3.1.4
is needed to establish tabularity, although it is not important if one is
only interested in the cell datum (see also Remark 4.2.5 below).


In \S4, we study examples of cellular algebras  
with
basis diagrams consisting of
non-intersecting curves that are inscribed in a $k$-box and labelled by
elements of a certain ring. 
We call the associated algebras ``decorated Temperley--Lieb algebras''.  

The original Temperley--Lieb algebras were defined by generators and relations,
together with key representations, in \cite{{\bf 34}}.  They are quotients of
the Hecke algebras associated to the symmetric groups.  The Temperley--Lieb
{\it diagram} algebra given in Definition 3.2.3 below 
is a realization of this
algebra, meaning that it is isomorphic to the algebra of \cite{{\bf 34}}.  We will
drop the word ``diagram'' for brevity.



\subhead 3.2 Temperley--Lieb diagrams \endsubhead

Note that in a Brauer diagram drawn on a single copy of the $k$-box it
is not generally possible to keep the embedded curves disjoint
(see Figure~2 for example). 
Let  $T_k(\emptyset) \subset B_k(\emptyset)$
denote the subset of diagrams having representative elements in which
the curves are disjoint. 
Representatives of this kind are called Temperley--Lieb diagrams.

It will be evident that $\Pl_n^B(\emptyset)$ has a subalgebra with
basis the subset $T_k(\emptyset) $. 
(That is to say, the disjointness property is preserved under multiplication.)
We denote this subalgebra 
$\Pl_n(\emptyset)$ 
(this may also be seen as 
a special case of \cite{{\bf 20}, Definition 1.8}).


Because of the disjointness property there is, for each element of 
 $T_k(\emptyset) $, a unique assignment of orientation to its 
curves that satisfies the following two conditions.

\item{(i)}{A curve meeting the $r$-th marked point of the standard
$k$-box, where $r$ is odd, must exit the box at that point.}

\item{(ii)}{Each connected component of the complement of the union of
the curves in the standard $k$-box may be oriented in such a way that
the orientation of a curve coincides with the orientation induced as
part of the boundary of the connected component.}

Note that the orientations match up automatically in composition. 

\definition{Definition 3.2.1} \cite{{\bf 7}}
Given a diagram, 
a point $x$ in the $k$-box is said to be $l$-exposed
(to the leftmost wall of the box) if $l$ is the smallest number of edges
it is necessary to traverse to get from $x$ to the leftmost wall.
Again because of disjointness, every point on an edge has the same
exposure. If this is $l$,
then the edge is said to be $l$-exposed. 

In composition each edge in the multiplied diagrams contributes a
segment (possibly all) to an edge in the product.  In this situation,
we call the edges in the multiplied diagrams the {\it ancestors} of
the corresponding edge in the product.  It will be evident
that the exposure of the new edge need not be the same as that of its
ancestors; however, the new exposure cannot exceed that of any
ancestor. 
\enddefinition
 


\ignore{{
We summarise the definition of the algebra $\Pl_k(\emptyset)$ from
\cite{{\bf 20}}.
}}

\ignore{{
\definition{Definition 3.1.2}
Let $k$ be a nonnegative integer.  Consider the set of 
finite collections of oriented disjoint curves 
(which we usually call ``edges''), 
smoothly embedded in the standard $k$-box,
with the property that there are precisely $k$ open curves (arcs) and
that the set of $2k$ endpoints coincides with the $2k$ marked points,
each endpoint meeting the boundary of the box transversely;
and that all curves are otherwise disjoint from the boundary of the box.
An element of $T_k(\emptyset)$,
called a diagram, consists of an equivalence class of
this set obtained after quotienting by smooth orientation-preserving 
diffeomorphisms of $\real^2$.
%
%
The orientations of the curves must satisfy the following two conditions.
\item{(i)}{A curve meeting the $r$-th marked point of the standard
$k$-box, where $r$ is odd, must exit the box at that point.}
\item{(ii)}{Each connected component of the complement of the union of
the curves in the standard $k$-box may be oriented in such a way that
the orientation of a curve coincides with the orientation induced as
part of the boundary of the connected component.}
\enddefinition
}}


\example{Example 3.2.2}
Let $k = 8$.  An element of $T_8(\emptyset)$ is shown in Figure~3.
Note that there are 10 connected components as in 
(ii) above, of which precisely 7 inherit a clockwise orientation.
\endexample

\topcaption{Figure~3} Typical element of $T_8(\emptyset)$ \endcaption
\centerline{
\hbox to 4.6in{
\vbox to 2.2in{\vfill
        \includegraphics{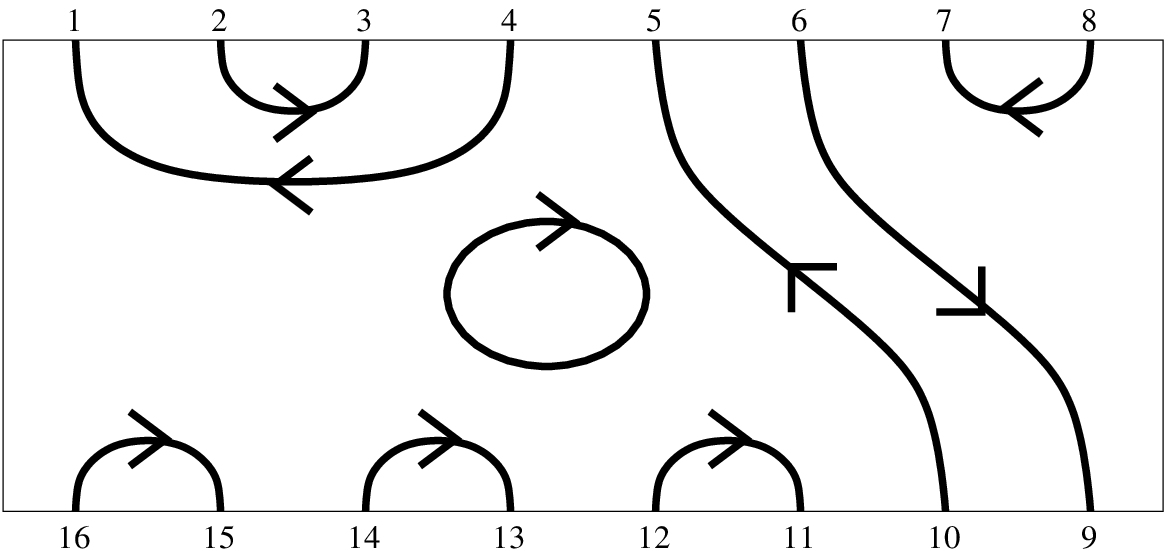}
}
\hfill}
}

\ignore{{
The following definition is a special case of \cite{{\bf 20}, Definition
1.8}.
\definition{Definition 3.1.4}
Let $k$ be a nonnegative integer.
The associative algebra $\Pl_k(\emptyset)$ over a commutative ring $R$
with $1$ is the free $R$-module having $T_k(\emptyset)$ as a basis,
with multiplication given by $k$-box multiplication. 
Note that the resulting curves are smooth.  The
orientations will match up automatically.
\enddefinition
}}



\definition{Definition 3.2.3}
Let $R$ be a commutative ring with $1$.  The Temperley--Lieb algebra,
$TL(n, \d)$, is the free $R[\d]$-module with basis given by the elements
of $T_n(\emptyset)$ with no closed loops.  The multiplication is
inherited from the multiplication on $\Pl_n(\emptyset)$ except that
one multiplies by a factor of $\d$ for each resulting closed loop and
then discards the loop.
\enddefinition

We usually consider $TL(n, \d)$ to be an algebra defined over $\A :=
\zed[v, v^{-1}]$, where $\d = v + v^{-1}$.  




\subhead 3.3 Decorated Temperley--Lieb algebras \endsubhead

We now recall from \cite{{\bf 20}, Example 2.2} the construction of
the algebra $P_n^A$ from the Temperley--Lieb algebra $TL(n, \d)$ and
the associative $R$-algebra $A$, where $R$ is a commutative ring
containing $\d$.  
The algebra $A$ is assumed to have identity and a trace functional 
$\tr: A \ra R$ with $\tr(ab) = \tr(ba)$ and $\tr(1) = \d$.

\definition{Definition 3.3.1}
Let $A$ be as above, and let $k$ be a nonnegative integer.
We define the tangles $T_k(A)$ to be those that arise from elements
of $T_k(\emptyset)$ by adding zero or more $1$-boxes labelled by
elements of $A$ to each edge.  
\enddefinition

Figure~4 shows a typical element of $T_8(A)$
in which $a, b, c, d, e \in A$.  

\topcaption{Figure~4} Typical element of $T_8(A)$ \endcaption
\centerline{
\hbox to 4.6in{
\vbox to 2.2in{\vfill
        \includegraphics{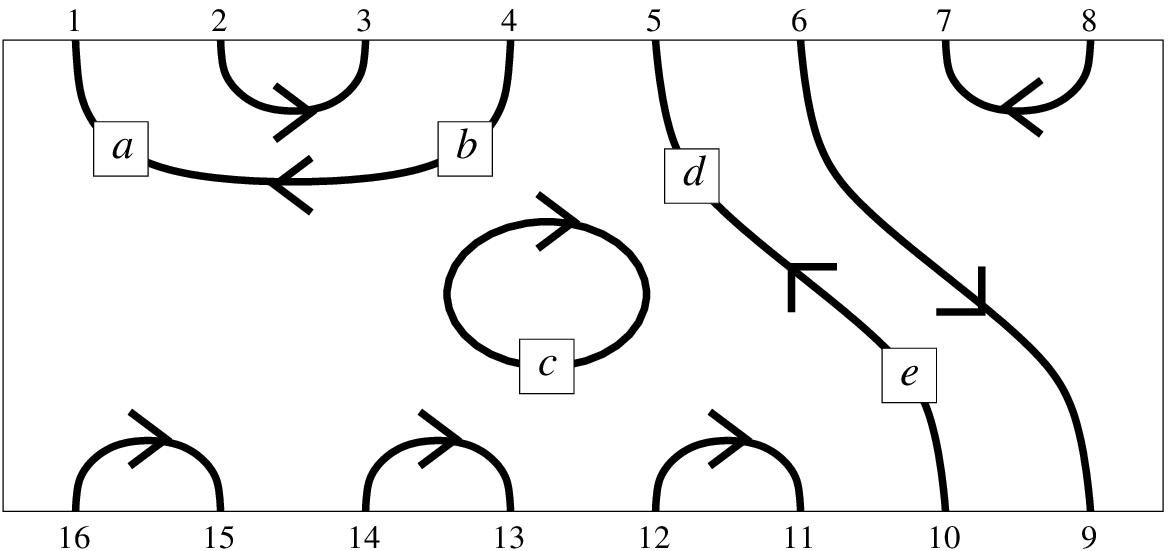}
}
\hfill}
}

Definition
3.2.3 generalizes naturally to this situation, as follows.

\definition{Definition 3.3.2}
Let $k$ be a nonnegative integer and let $A$ be an $R$-algebra (as
before) with a free $R$-basis, $B_A = \{a_i : i \in I\}$, where 
$1 \in \{a_i\}$.
The associative $R$-algebra $P_k^A$
is the free $R$-module having as a basis those elements of $T_k(A)$ 
satisfying the conditions that
\item{(i)}{all labels on edges are basis elements $a_i$,}
\item{(ii)}{each edge has precisely one label and}
\item{(iii)}{there are no closed loops}.

The multiplication is defined   
on basis elements of $P_k^A$ as above, and extended bilinearly.
We start with the multiplication on $T_k(A)$, then 
apply
relations (a), (b) and (c) below to express the product as an
$R$-linear combination of basis elements, and finally, apply relation (d)
below to remove any loops, multiplying by the scalar shown for each
loop removed.
\enddefinition

\topcaption{Figure~5} Relation (a) of Definition 3.3.2 \endcaption
\centerline{
\hbox to 3.3in{
\vbox to 0.9in{\vfill
        \includegraphics{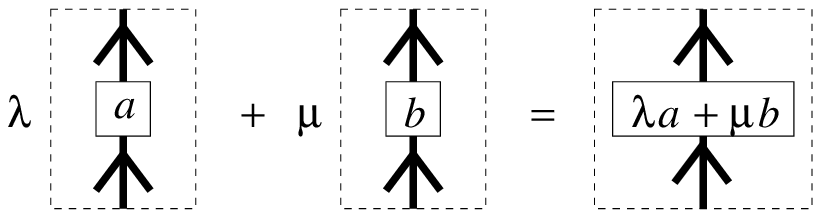}
}
\hfill}
}

\topcaption{Figure~6} Relation (b) of Definition 3.3.2 \endcaption
\centerline{
\hbox to 1.5in{
\vbox to 1.4in{\vfill
        \includegraphics{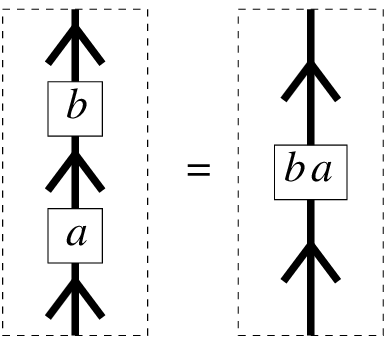}
}
\hfill}
}
\centerline{
\hbox to 1.5in{
\vbox to 1.4in{\vfill
        \includegraphics{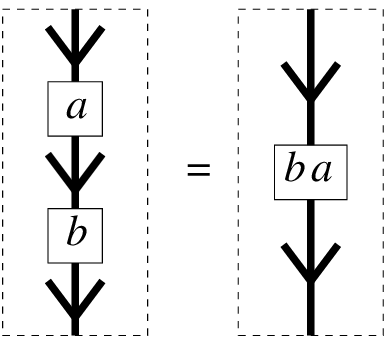}
}
\hfill}
}

\topcaption{Figure~7} Relation (c) of Definition 3.3.2 \endcaption
\centerline{
\hbox to 1.5in{
\vbox to 1.4in{\vfill
        \includegraphics{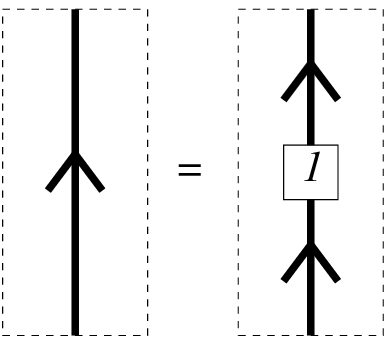}
}
\hfill}
}


\topcaption{Figure~8} Relation (d) of Definition 3.3.2 \endcaption
\centerline{
\hbox to 1.7in{
\vbox to 1.3in{\vfill
        \includegraphics{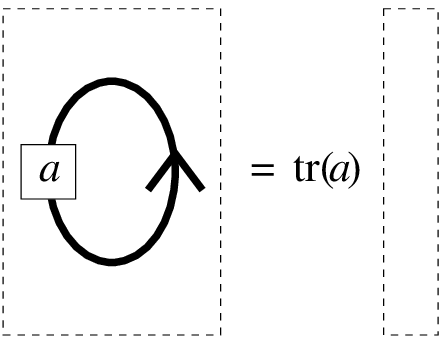}
}
\hfill}
}

We call the algebra $P_k^A$ a {\it decorated Temperley--Lieb algebra} and
the above basis, denoted $\BB_n^A$,  
the {\it canonical basis} with respect to $B_A$.

For a proof that this procedure does define an associative algebra, the 
reader is referred to \cite{{\bf 20}, Example 2.2}. 

\remark{Remark 3.3.3}
The direction on the arrow in relation (d) is immaterial, although one
can define a more intricate version of the algebra in which there are
two traces on $A$, one for each orientation of the arrow.
\endremark

As mentioned in 
\cite{{\bf 20}}, this construction might be regarded as a kind of wreath
product of $TL(n)$ with $A$.


Fix a natural number $l$. 
Consider the subset of the canonical basis consisting of elements 
with the property that every edge with exposure greater than $l-1$ is
decorated by the 1-box containing the identity element of $A$. 
It follows from the definition of exposure that this property is
preserved under multiplication, and hence that the subset, denoted
$\BB_n^{A,l}$,
forms a basis for a subalgebra $P_n^{A,l}$ of $P_n^A$. 

For example, $P_n^{A,0}$ is isomorphic to the ordinary Temperley--Lieb
algebra, while $P_n^{A,n} \cong P_n^A$. 
The algebras $P_n^{A,l}$ come from the contour algebra formalism
introduced in \cite{{\bf 7}}. 

\proclaim{Proposition 3.3.4}
Let $A$ be a hypergroup over $R$ with distinguished
basis $B$ and any trace map.  There is an isomorphism $\rho$ of $R$-algebras
from $(A^{\otimes n}, B^{\otimes n})$ to the subalgebra of $P_n^A$ spanned by
all canonical basis elements with no non-propagating edges.  The
isomorphism takes basis elements to basis elements.
\endproclaim

\demo{Proof}
Let $b = b_{i_1} \otimes b_{i_2} \otimes \cdots \otimes b_{i_n}$ be a
typical basis element from the set $B^{\otimes n}$.  This element is
sent by the isomorphism, $\rho$, to a canonical basis element of $P_n^A$ with
no non-propagating edges, where the decoration on the $k$-th propagating
edge (counting from $1$ to $n$, starting at the left) is $b_{i_k}$ if
$k$ is odd, and $b_{\overline{i_k}}$ if $k$ is even.  

For the proof that this construction defines an isomorphism of algebras,
the reader is referred to \cite{{\bf 15}, Proposition 2.3.4}.
Note that the trace map plays no role in the structure of the algebra,
because closed loops cannot arise.
\qed\enddemo

\proclaim{Lemma 3.3.5 \cite{{\bf 15}, Lemma 2.3.2}}
Let $A$ be a hypergroup over $\zed$ with distinguished
basis $B$ and the trace map $\d . t$, where $t$ is as in Proposition 1.2.4.
There is an linear anti-automorphism, $*$, of $P_n^A$ permuting the canonical
basis.  The image, $b^*$, of a basis element $b$ under this map is
obtained by reflecting $b$ in the line $y = 1/2$, reversing the
direction of all the arrows and replacing each $1$-box labelled by
$b_i \in B$ by a $1$-box labelled by $b_{\bar{i}}$.
\qed\endproclaim


\head 4. Tabularity of decorated Temperley--Lieb algebras and related
algebras \endhead

\subhead 4.1 General results \endsubhead

The following result is \cite{{\bf 15}, Theorem 3.2.3}; we recall part of the proof
below as we require the construction later.

\proclaim{Theorem 4.1.1 \cite{{\bf 15}}}
Let $A$ be a hypergroup over $\zed$ with distinguished
basis $B$ and the trace map $\d . t$, where $t$ is as in Proposition 1.2.4.
Then the algebra $P_n^A$ equipped with its canonical basis $\BB_n^A$ 
is a tabular algebra.
\endproclaim

\demo{Proof}
We require to construct a table datum. 
Let $\Lambda$ be the set of integers $r$ with $0
\leq r \leq n$ and $n - r$ even, ordered in the usual way.  

For $\l \in \Lambda$, let $(\G(\l), B(\l))$ be the $\l$-th tensor
power of the hypergroup $(A, B)$ with the basis and
anti-automorphism induced by Proposition 1.2.5.

Let $M(\l)$ be the set of possible configurations of $(n-\l)/2$ non-propagating
edges with endpoints on the line $y = 1$ that arise from an element of 
$\BB_n^A$.  Let $b = b_{i_1} \otimes b_{i_2} \otimes \cdots \otimes
b_{i_\l}$ be a
basis element of $B(\l)$ and let $m$ and $m'$ be elements of $M(\l)$.
The map $C$ produces a basis element in $\BB_n^A$ from the triple 
$(m, b, m')$ as follows.  Turn the half-diagram corresponding to $m'$
upside down, reverse the directions of all the arrows and relabel all
$1$-boxes labelled by $b_i \in B$ so they are labelled by
$b_{\bar{i}}$.
Join any free marked points in the line $y = 0$ 
to free marked points in 
the line $y = 1$ so that they do not intersect.  Orient any new edges
according to the orientation of the standard $n$-box.  Decorate the $\l$
propagating edges with the basis element $b$ using the construction
of Proposition 3.3.4.  (See \cite{{\bf 15}, Example 3.2.4} for an illustration
of this.)

The map $*$ is the one given by Lemma 3.3.5.

For the proof that this construction defines a table datum, the reader is
referred to \cite{{\bf 15}, Theorem 3.2.3}, which proves the stronger result
that the construction defines a tabular algebra with trace.
\qed\enddemo


Similarly, using the subdatum idea we have
\proclaim{Theorem 4.1.2}
Let $A$ be a hypergroup over $\zed$ with distinguished
basis $B$ and the trace map $\d . t$, where $t$ is as in Proposition 1.2.4.
Then the algebra $P_n^{A,l}$ equipped with its canonical basis $\BB_n^{A,l}$ 
is a tabular algebra (for any appropriate $l$).
\endproclaim

\demo{Proof}
The proof goes through as before (noting that the hypergroups
appearing in the table datum are the $\lambda$-th tensor power of
$(A, B)$ for $\lambda < l$, and the $l$-th tensor power thereafter). 
\qed\enddemo

\proclaim{Theorem 4.1.3}
Suppose that $R$ is a commutative ring with identity and $(A, B)$ 
is a hypergroup such that $R \otimes_\zed A$ is cellular with respect 
to an anti-automorphism of the based ring $A$.  Equip $A$ with the 
trace map $\d . t$, where $t$ is the trace of Proposition 1.2.4.  
Then $R \otimes_\zed P_n^A$ is cellular, 
and an explicit cell datum is given by Theorem 2.2.1.
\endproclaim

\demo{Proof}
By Proposition 1.3.5, the anti-isomorphism of the statement is of the form
$\bar{g}$ for some automorphism $g$ of based rings.  The map $g$ then induces
a permutation $\a$ of the basis diagrams of $P_n^A$ by acting simultaneously on
the decorations of each edge; note that $\a$ preserves the number of
propagating edges in each diagram.  By Proposition 1.2.4, we see that 
$t(g(x)) = t(x)$, because the action of $g$ on an element $x \in A$ does
not alter the coefficient of the identity element of $\BB$.  The relations
of Definition 3.3.2 now show that $\a$ induces a basis-preserving 
automorphism of the algebra $R \otimes_\zed P_n^A$.

As explained in Theorem 4.1.1, the hypergroups occurring in the table
datum for $P_n^A$ are certain tensor powers of the hypergroup $(A, \BB)$ (see
Proposition 1.2.5).  The automorphism $\a$ induces the basis-preserving
automorphism $g^r$ on the hypergroup $(A^{\otimes r}, \BB^r)$, and the
latter hypergroup is cellular over $R$ with respect to the anti-automorphism
$\overline{g^r}$ by Example 1.1.6.

The result now follows by Theorem 2.2.1, in which $\a$ is as above and 
$f_r$ is given by $g^r$.
\qed\enddemo


\subhead 4.2 Cyclotomic Temperley--Lieb algebras \endsubhead

The so-called cyclotomic Temperley--Lieb algebras $TL_{n, m}(\d_0, \d_1, 
\ldots\d_{m-1})$ are algebras over a ring $R$ containing elements
$\d_0, \d_1, \ldots, \d_{m-1}$.  They were introduced by 
Rui and Xi \cite{{\bf 31}}, both in terms of generators 
and relations \cite{{\bf 31}, Definition 2.1}, and
equivalently 
in terms of a calculus of diagrams \cite{{\bf 31}, Definition~3.3, 
Theorem~3.4}. 

The rules for manipulating the decorations on edges in \cite{{\bf 31}} are somewhat
intricate (see \cite{{\bf 31}, \S3}). 
We show here, for certain values
of the parameters, how the algebra may also 
be set up using planar algebras on $1$-boxes.  

\remark{Remark 4.2.1} This construction is mentioned
in passing by Cox et al \cite{{\bf 7}, Remark 2.3}, who then provide an entirely
straightforward construction for an isomorphic diagram algebra and a sequence
of generalizations, which we will not need here.
It is also possible to define
the algebras for general parameter values using planar algebras on $1$-boxes,
but this requires the use of two traces (see Remark 3.3.3) and we will not 
pursue this here.
\endremark

\proclaim{Lemma 4.2.2}
Let $(A, B)$ be the cyclic group $\zed_m$ considered as a hypergroup 
$(\zed \zed_m, \zed_m)$, equipped with the trace map $\d . t$, where $t$
is as in Proposition 1.2.4.  Then $P_n^A$ is isomorphic to the cyclotomic
Temperley--Lieb algebra $TL_{n, m}(\d, 0, 0, \ldots, 0)$ over $\A$, and the
canonical basis of $P_n^A$ can be chosen to map to the diagram basis
of $TL_{n, m}$.
\endproclaim

\demo{Proof}
For $1 \leq k < n$, let $E_{k, n}$ be an element of $T_n(\emptyset)$ 
with no closed loops
in which each  point $i$ is connected by a propagating edge to point 
$2n + 1 - i$, unless $i \in \{k, k+1, 2n - k, 2n + 1 - k\}$; furthermore, 
points $k$ and $k+1$ are connected by a non-propagating edge, as are 
points $2n - k$ and $2n + 1 - k$.

Let $g$ be a generator of the group $\zed_m$.  
For $1 \leq k \leq n$, let $T_{k, n}$ be the basis diagram corresponding to $$
\underbrace{1 \otimes 1 \otimes \cdots \otimes 1}_{k-1} 
\otimes g \otimes
\underbrace{1 \otimes 1 \otimes \cdots \otimes 1}_{n-k} 
$$ under the isomorphism of Proposition 3.3.4.

The map sending $E_{i, n}$ and $T_{j, n}$ to the respective elements
$E_i$ and $T_j$ in $TL_{n, m}$
(see the notation of the proof of \cite{{\bf 31}, Theorem 3.4}) extends uniquely
to an isomorphism of $\A$-algebras.  This is a matter of checking that
the multiplications in the two diagram calculi are compatible, and this 
follows from the rules given in \cite{{\bf 31}, \S3}.
\qed\enddemo

The remarks in the proof of Theorem 4.1.1 now give

\proclaim{Corollary 4.2.3}
The cyclotomic Tem\-per\-ley--Lieb algebra $TL_{n, m}(\delta, 0, \ldots, 0)$
is a tabular algebra with trace.  
\qed\endproclaim

In contrast, for general parameter values, the cyclotomic Temperley--Lieb
algebra is not tabular in any obvious way.  This is not so surprising
given that it is a multiparameter algebra, but even if we require all
the parameters to lie in $\A$, complications arise.
If $\d_i = \d_{m-i} \in \A$ for all $0 \leq i \leq m$, 
there is an $\A$-linear anti-automorphism of the algebra fixing the generators
$e_i$ and sending each $t_j$ to $t_j^{-1}$, and this serves as a tabular 
anti-automorphism.  In general, however, this map fails to be an $\A$-linear 
anti-automorphism, even if all the $\d_i$ lie in $\A$.  If one is primarily
interested in the cellular structure, this is not a major problem, as we
will explain in Remark 4.2.5.

\proclaim{Corollary 4.2.4 (Rui--Xi, \cite{{\bf 31}, Theorem 5.3})}
Let $R$ be a commutative ring with identity such that $x^m - 1$ splits into
linear factors over $R[x]$.  Then the cyclotomic Temperley--Lieb algebra
$TL_{n, m}(\d, 0, 0, \ldots, 0)$ over $R$ is cellular with respect to the
map $* \circ \a = \a \circ *$, where $*$ is as defined in
Lemma 3.3.5, and $\a$ is the map induced by applying the inversion
automorphism of $\zed_m$ to each edge of each basis diagram.
\endproclaim

\demo{Proof}
Let $P(x) = x^m - 1 \in R[x]$.  By hypothesis, $P(x)$ splits into linear
factors $$P(x) = \prod_{i = 1}^k (x - \l_i).$$  By Example 1.1.4, this shows 
that $R \otimes_\zed \zed\zed_m = R[x]/\lan P(x) \ran$ is cellular with 
respect to the identity map.  Since the identity map $\zed_m \ra \zed_m$ is an
anti-automorphism of based rings, Theorem 4.1.3 and Lemma 4.2.2 show that
$TL_{n, m}(\d, 0, 0, \ldots, 0)$ is cellular, and Theorem 2.2.1 provides a
cell datum with the required properties.
\qed\enddemo

\remark{Remark 4.2.5}
Rui and Xi \cite{{\bf 31}, Theorem 5.3} prove the result above for arbitrary 
values of the parameters $\d_i$, but the cell datum remains essentially 
the same in each 
case.  In particular, the map $* \circ \a = \a \circ *$ 
remains as the cellular anti-automorphism of $TL_{n, m}$ for all parameter 
values, even though $*$ (respectively, $\a$) is not an anti-automorphism 
(respectively, an automorphism) of the algebra in general.
\endremark


\subhead 4.3 Other similar examples, and subdata \endsubhead

Other algebras that can be treated similarly 
include the blob algebra and the
generalized Temperley--Lieb algebra of type $H_n$.  
We will give only a sketch of the arguments, for the sake of brevity. 

The blob algebra was defined by the second author and Saleur \cite{{\bf 30}} 
in a statistical mechanical context.  
It may be defined as a certain subalgebra of
$P_n^A$, where $A$ is the algebra $\zed[x]/\lan x^2 - x \ran$.  Since $A$
has no obvious hypergroup structure, we make the change of variables
$y = 2x - 1$ and consider $A' = \zed[y]\lan y^2 - 1 \ran$; this is
the hypergroup $(\zed \zed_2, \zed_2)$.  (Note that
$\kyu \otimes_\zed A \cong \kyu \otimes_\zed A'$; this has the effect of
making the decoration in \cite{{\bf 30}} unipotent instead of idempotent.)
The blob algebra can be constructed from $P_n^{A'}$ using a 
subdatum to ensure
that only edges exposed to the leftmost edge of the $n$-box may carry
decorations; this means that most of the hypergroups appearing in the table
datum are isomorphic to $A'$.  Since $A'$ is cellular with respect to 
the identity map (\idest inversion in $\zed_2$), Theorem 2.2.1 can be
applied with $\a$ taken to be the identity automorphism.

The generalized Temperley--Lieb algebras of type $H_n$ were considered
by Graham \cite{{\bf 10}} and an explicit cell datum was constructed
in \cite{{\bf 13}, Theorem 3.3.5} in terms of a basis of diagrams that were
later shown to be a tabular basis \cite{{\bf 14}, Theorem 5.2.5}.
The treatment of these algebras using Theorem 2.2.1 is similar to that 
of the blob algebra above (see also the remarks in \cite{{\bf 14}, \S2.1}).
The main differences are 
(a) the relevant hypergroup to use is $A = \Gamma_H = 
\zed[x]/\lan x^2 - x - 1 \ran$ with basis $B_H$ given by the images of 
$1$ and $x$,
and 
(b) more care is needed in defining the subdatum, which may 
be constructed using the rules for ``H-admissibility'' listed in 
\cite{{\bf 15}, Definition 4.2.3}.

\head 5. Tabularity of the Brauer algebra and related algebras \endhead

\subhead 5.1 The Brauer algebra \endsubhead

We now recall the tabular structure of the Brauer algebra.
This example comes from \cite{{\bf 14}, Example 2.1.2} and \cite{{\bf 16}, \S4.2}.


As in \cite{{\bf 11}, \S4}, we may describe the basis diagrams in terms of
certain triples.  

\definition{Definition 5.1.1}
Fix a Brauer diagram $D$. 
The integer $t(D)$ is defined to be the number of propagating edges.
The involutions $S_1(D)$, $S_2(D)$ in the symmetric
group $\Sy(n)$ are defined such that $S_i(D)$ interchanges the ends of
the joins between points with the same $y$-coordinate.  For example, with
$D$ as in Figure~2, we have $S_1(D) = (1 6)(2 5)$.
Corresponding to these we have subsets $\text{\rm Fix}(S_i(D))$ of $\{1,
\ldots, n\}$, which are the fixed points of the involutions $S_i(D)$.
Finally, we have a permutation $w(D)$ in $\Sy(t)$,
where $t = t(D)$; this is the permutation of $\text{\rm Fix}(S_1(D))$
determined by taking the end points of the propagating edges (regarded
as joining from the row $y = 0$ to the row $y = 1$) in the order determined 
by taking their starting points in the row $y = 0$ in increasing order.  
(We consider $\Sy(0)$ to be the trivial group, in which case $w$ is the 
identity.)  The diagram $D$ is then determined by the triple 
$[S_1(D), S_2(D), w(D)]$.
\enddefinition


A table datum for the Brauer algebra (equipped with the diagram basis)
may be given as follows.  (This gives the algebra the structure of a
tabular algebra with trace.)

\definition{Definition 5.1.2}
Let $B_n$ be the Brauer algebra (over $\A$) on $n$ strings.  
The algebra has a table datum $(\Lambda, \Gamma, B, M, C, *)$ as follows.

Take $\Lambda$ to be the set of integers $i$ between $0$ and $n$ such that
$n - i$ is even, ordered in the natural way.  
If $\l = 0$, take $(\Gamma(\l), B(\l))$ to be the
trivial hypergroup; otherwise, take $\Gamma(\l)$ to
be the group ring $\zed \Sy(\l)$ with basis $B(\l) = \Sy(\l)$ and 
involution $\overline{w}
= w^{-1}$.  Take $M(\l)$ to be the set of involutions on $n$ letters 
with $\l$ fixed points.  Take $C(S_1, w, S_2) = [S_1,
S_2, w]$; $\Im(C)$ contains the identity element.  
The anti-automorphism $*$ sends $[S_1, S_2, w]$ to $[S_2, S_1, w^{-1}]$.
\enddefinition

We state the next result for later use.

\proclaim{Lemma 5.1.3}
\item{\rm (i)}{The operation of reflecting each basis diagram of $B_n$ in 
a vertical line $x = (n+1)/2$ extends to a unique automorphism of 
$\A$-algebras, $\rho$, of $B_n$.}
\item{\rm (ii)}{Let $\omega_k : \Sy_k \ra \Sy_k$ be the automorphism of the
symmetric group obtained by conjugation by the longest element of $\Sy_k$.
Then we have $$\rho([S_1, S_2, w]) = 
[\omega_n(S_1), \omega_n(S_2), \omega_t(w)],
$$ where $t$ is the number of fixed points of $S_1$ and $S_2$.}
\endproclaim

\demo{Proof}
Part (i) follows easily from the definition of the multiplication in $B_n$.

Part (ii) is a consequence of the observation that if $g \in \Sy_k$, we
have $g(i) = j$ if and only if $(\omega_k(g))(n + 1 - i) = n + 1 - j$.
\qed\enddemo

\remark{Remark 5.1.4}
As mentioned in \cite{{\bf 14}, \S2.1}, similar techniques
may be applied to the case of the
partition algebra of \cite{{\bf 28}}; again the hypergroups are
symmetric groups equipped with inversion as the involution.  This
recovers Xi's main result in \cite{{\bf 36}}.  
\endremark


\subhead 5.2 The walled Brauer algebra \endsubhead

The walled Brauer algebra, also known as the rational Brauer algebra, is
a certain subalgebra of the Brauer algebra first considered by Benkart et
al in \cite{{\bf 5}, \S5}.  The cellularity of this algebra is well-known to the
experts, and it is implicit in the construction of the basis described in 
\cite{{\bf 23}}.  We include the example here to illustrate how easy it is to
describe the cellular structure of this algebra using our techniques.

\definition{Definition 5.2.1}
Let $(\Lambda, \Gamma, B, M, C, *)$ be the table datum for the Brauer
algebra $B_n$ given in Definition 5.1.2, let $r$ and $s$ be positive integers
with $r + s = n$, and define $(\Lambda', \Gamma', B', M', C', *')$, 
as follows.  Let $$
\Lambda' = \{\l' \in \Lambda : \l' \geq |r - s|\}
,$$ and for each $\l' \in \Lambda'$, define $(\Gamma'(\l'), B'(\l'))$
to be the hypergroup corresponding to the subgroup $\Sy(r') \times \Sy(s')$
of $\Sy(\l')$, where $r'$ and $s'$ are the unique nonnegative integers 
satisfying $r' + s' = \l'$ and $r' - s' = r - s$.  For $\l' \in
\Lambda'$ and the corresponding integers $r'$ and $s'$ just given, denote
by $P^-$ the set $\{1, 2, \ldots, r\}$ and by $P^+$ the set
$\{r+1, r+2, \ldots, n\}$.  We define $M'(\l')$ to be the subset of 
$M(\l')$ consisting of involutions $S$ for which the following conditions
hold: 
\item{(i)}{$S$ has $r'$ fixed points in $P^-$;}
\item{(ii)}{$S$ has $s'$ fixed points in $P^+$;}
\item{(iii)}{if $S$ exchanges two distinct points then one of the points
comes from $P^-$ and the other from $P^+$.}

We define $C'$ and $*'$ to be the restrictions of $C$ and $*$ to the
appropriate domains.
\enddefinition

\proclaim{Definition/Lemma 5.2.2}
The tuple $(\Lambda', \Gamma', B', M', C', *')$ is a subdatum for the table
datum of the Brauer algebra, and corresponds to an algebra $A'$, which is
by definition the {\it walled Brauer algebra}, $B_{r, s}$.
\endproclaim

\demo{Proof}
The walled Brauer algebra $B_{r, s}$ defined in \cite{{\bf 5}, \S5} is 
given to be that spanned by certain basis diagrams, called $(r, s)$-diagrams.
A routine check shows that the $(r, s)$-diagrams are precisely those
in the image of the map $C'$.  It remains to be checked that $\Im(C')$ is
a subalgebra of $A$, but this also presents no difficulties (see
\cite{{\bf 5}, \S5}).
\qed\enddemo

\proclaim{Corollary 5.2.3}
The walled Brauer algebra $B_{r, s}$ is cellular with respect to the 
anti-automorphism $*'$ of Lemma 5.2.2.
\endproclaim

\demo{Proof}
Since $B_{r, s}$ contains the identity element of the Brauer algebra
$B_{r+s}$, Proposition 2.1.5 and Lemma 5.2.2 show that 
$(\Lambda', \Gamma', B', M', C', *')$ is a table datum for $B_{r, s}$
(and furthermore, that $B_{r, s}$ is a tabular algebra with trace).

Example 1.1.3 shows that $\zed \Sy_n$ is cellular over $\zed$ with
respect to the group inversion.  Example 1.1.6 then shows that
$\zed \Sy_r \otimes_\zed \zed \Sy_s \cong \zed (\Sy_r \times \Sy_s)$
is cellular with respect to inversion in the group $\Sy_r \times \Sy_s$.

Using the above observations, Theorem 2.2.1 (with $\a$ as the identity map)
now constructs a cell datum showing that $B_{r, s}$ is cellular with 
respect to $*' \circ \a = *'$, as required.
\qed\enddemo

\remark{Remark 5.2.4}
Notice that Theorem 5.2.3 includes the cellularity of the Brauer algebra
as a special case; this was originally a theorem of Graham--Lehrer
\cite{{\bf 11}, Theorem 4.10}.
\endremark


\subhead 5.3 Jones' annular algebra \endsubhead

Jones' annular algebra (or the Jones algebra, for short) 
is a certain subalgebra of the Brauer algebra that
is also a quotient of an affine Hecke algebra of type $A$.  It was
introduced in \cite{{\bf 19}} and first shown to be cellular in \cite{{\bf 11}, \S6}.
In order to define the algebra, we recall the notion of an annular
involution.

\definition{Definition 5.3.1}
An involution $S \in \Sy_n$ is annular if and only if for each pair $i, j$
interchanged by $S$ (with $(i < j)$) and $P_{i, j} = \{k : i \leq k \leq j\}$, we have
\item{(a)}{$S(P_{i, j}) = P_{i, j}$ and}
\item{(b)}{either $S$ fixes no element of $P_{i, j}$ or every element fixed
by $S$ is contained in $P_{i, j}$.}
\enddefinition

\definition{Definition 5.3.2}
Let $(\Lambda, \Gamma, B, M, C, *)$ be the table datum for the Brauer algebra
given by Definition 5.1.2.  We define $\Lambda' = \Lambda$, and for each
$\l' \in \Lambda$, let $(\Gamma(\l'), B(\l'))$ be the cyclic group of order
$\l'$ considered as a hypergroup, unless $\l' = 0$, in which case we define
 $(\Gamma(\l'), B(\l'))$ to be the trivial hypergroup.  Also for each 
$\l' \in \Lambda'$, we define $M(\l')$ to be the set of annular involutions
with $\l'$ fixed points; it may be checked that this is always a nonempty set.
We define $C'$ and $*'$ to be the restrictions of $C$ and $*$ to the
appropriate domains.
\enddefinition

\proclaim{Definition/Lemma 5.3.3}
The tuple $(\Lambda', \Gamma', B', M', C', *')$ is a subdatum for the table
datum of the Brauer algebra, and corresponds to an algebra $A'$, which is
by definition Jones' annular algebra, $J_n$.
\endproclaim

\demo{Proof}
The definition of $J_n$ given in \cite{{\bf 11}, \S6} is as the span of those
Brauer algebra diagrams $[S_1, S_2, w]$, where $S_1$ and $S_2$ are annular
involutions with $t$ fixed points, and where $w$ is an element of the cyclic
group of order $t$ if $t > 0$, with $w = 1$ if $t = 0$.  It is not hard to 
see that this agrees with our construction.  It is also routine to check
that this defines a subalgebra; see \cite{{\bf 19}} or \cite{{\bf 11}, \S6}.
\qed\enddemo

Note that the tabular structure of the Jones algebra has already been
described in \cite{{\bf 14}, Example 2.1.4}.

\example{Example 5.3.4}
Let $[S_1, S_2, w]$ be the basis diagram shown in Figure~2.  The involution
$S_1$ is annular because the subsets $\{k : 2 \leq k \leq 5\}$ and
$\{k : 1 \leq k \leq 6\}$ contain all the fixed points of $S_1$.  The
involution $S_2$ is annular because $S_2$ fixes no elements in the sets
 $\{k : 2 \leq k \leq 3\}$ and $\{k : 5 \leq k \leq 6\}$.  Note that
$w$ lies in the cyclic group of order $2$, the number of fixed points of
each of $S_1$ and $S_2$.
\endexample

The reason for the term ``annular'' is that the basis diagrams of $J_n$ are
precisely those that can be inscribed without intersections within an annulus.
Although the diagram in Figure~2 has intersections when inscribed in a
rectangle, it can be inscribed without intersections in an annulus, as shown
in Figure~9.

\topcaption{Figure~9} A basis element of Jones' annular algebra \endcaption
\centerline{
\hbox to 2.300in{
\vbox to 2.000in{\vfill
        \includegraphics{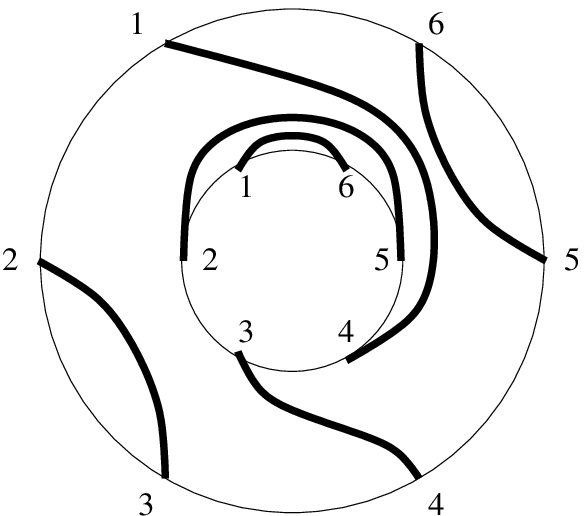}
}
\hfill}
}

A cell datum for the Jones algebra cannot be obtained by restriction of
the cell datum for the Brauer algebra, the obstruction being essentially
that group algebras of cyclic groups are not usually cellular with respect 
to the group inversion.  However, we can use our main results to exploit 
the fact that these group algebras are cellular with respect to the identity
anti-automorphism (see Example 1.1.4 and the proof of Corollary 4.2.4).

\proclaim{Lemma 5.3.5}
The automorphism $\rho$ of the Brauer algebra $B_n$, defined in Lemma 5.1.3,
restricts to an automorphism of $J_n$.
\endproclaim

\demo{Proof}
If $S_1$ is an annular involution, it follows from the symmetric nature of
Definition 5.3.1 that $\omega_n(S_1)$ is also annular.  If $w \in \zed_t$,
we find that $\omega_t(w) = w^{-1}$.  The result now follows from 
Lemma 5.1.3 (ii).
\qed\enddemo

\proclaim{Corollary 5.3.6 (Graham--Lehrer, \cite{{\bf 11}, Theorem 6.15})}
Let $n \in \enn$, and let $R$ be a commutative ring with identity such 
that $x^t - 1$ splits into linear factors over $R[x]$ for all 
$0 \leq t \leq n$ such that $n - t$ is even.
Then the Jones algebra $J_n$ over $R[v, v^{-1}]$ is cellular with respect 
to the map $* \circ \rho = \rho \circ *$, where $*$ is the map of 
Definition 5.1.2, and $\rho$ is the map of Lemma 5.3.5.
\endproclaim

\demo{Proof}
Using Example 1.1.4 as in the proof of Corollary 4.2.4, we see that
$R \zed_t$ is cellular over $R$ with respect to the identity map, for all
values of $t$ given in the statement.  We then apply Theorem 2.2.1 with
$\rho$ in the role of $\a$.  Note that in this case, the maps $f_\l$ all
arise from inversion in suitable cyclic groups, and the maps $\s_\l$ are
the maps $\omega_n$ from Lemma 5.1.3.  The composite maps $\overline{f_\l}$
are all equal to identity maps.  Theorem 2.2.1 completes the proof by
constructing a cell datum.
\qed\enddemo

\remark{Remark 5.3.7}
It is possible to exploit the rotational symmetry in the Jones algebra
to define automorphisms of $J_n$ other than $\rho$ for which Corollary
5.3.6 still holds.  This shows that the automorphism $\a$ required by
Theorem 2.2.1 need not be uniquely determined.
\endremark

\remark{Remark 5.3.8}
Our approach in this paper provides a convenient framework for examining
situations such as the embedding of the Jones algebra in the Brauer algebra,
where the cellular structures are not compatible, but the tabular structures
are.
\endremark




\head Acknowledgements \endhead

The first author thanks Stephen Doty and Steffen Koenig for helpful 
conversations and correspondence.  The authors 
also thank the referee for helpful comments.



\leftheadtext{} \rightheadtext{}

\end
\Refs\refstyle{A}\widestnumber\key{[30]} \leftheadtext{References}
\rightheadtext{References}

\ref\key{{\bf 1}}
\by Z. Arad and H.I. Blau
\paper On Table Algebras and Applications to Finite Group Theory
\jour J. Algebra
\vol 138 \yr 1991 \pages 137--185
\endref

\ref\key{{\bf 2}}
\by Z. Arad, E. Fisman and M. Muzychuk
\paper Generalized table algebras
\jour Isr. J. Math.
\vol 114 \yr 1999 \pages 29--60
\endref

\ref\key{{\bf 3}}
\by S. Ariki
\paper Robinson--Schensted correspondence and left cells
\jour Adv. Studies Pure Math.
\vol 28 \yr 2000 \pages 1--20
\endref

\ref\key{{\bf 4}}
\by E. Bannai and T. Ito
\book Algebraic Combinatorics I: Association Schemes
\publ Ben\-jam\-in--Cum\-mings
\publaddr Menlo Park, CA
\yr 1984
\endref

\ref\key{{\bf 5}}
\by G. Benkart, M. Chakrabarti, T. Halverson, R. Leduc, C. Lee and J. Stroomer
\paper Tensor product representations of general linear groups and their 
connections with Brauer algebras
\jour J. Algebra \vol 166 \yr 1994 \pages 529--567
\endref

\ref\key{{\bf 6}}
\by R. Brauer
\paper On algebras which are connected with the semisimple continuous groups
\jour Ann. of Math.
\vol 38
\yr 1937
\pages 854--887
\endref

\ref\key{{\bf 7}}
\by A. Cox, P. Martin, A. Parker and C. Xi
\paper Representation theory of towers of recollement: theory, notes and
examples
\miscnote preprint; {\tt math.RT/0411395}
\endref

\ref\key{{\bf 8}}
\by J. Du
\paper IC Bases and Quantum Linear Groups
\jour Proc. Sympos. Pure Math.
\vol 56 \yr 1994 \pages 135--148
\endref

\ref\key{{\bf 9}}
\by J. Enyang
\paper Cellular bases for the Brauer and Birman-­Murakami-­Wenzl algebras
\jour J. Algebra 
\vol 281 \yr 2004 \pages 413-­449
\endref

\ref\key{{\bf 10}}
\by J.J. Graham
\book Modular representations of Hecke algebras and related algebras
\publ Ph.D. thesis
\publaddr University of Sydney
\yr 1995
\endref

\ref\key{{\bf 11}}
\by J.J. Graham and G.I. Lehrer
\paper Cellular algebras
\jour Invent. Math.
\vol 123
\yr 1996
\pages 1--34
\endref

\ref\key{{\bf 12}}
\by R.M. Green
\paper Completions of cellular algebras
\jour Comm. Alg.
\vol 27 \yr 1999 \pages 5349--5366
\endref

\ref\key{{\bf 13}}
\by R.M. Green
\paper Cellular algebras arising from Hecke algebras of type $H_n$
\jour Math. Zeit.
\vol 229 \yr 1998 \pages 365--383
\endref

\ref\key{{\bf 14}}
\by R.M. Green
\paper Tabular algebras and their asymptotic versions
\jour J. Algebra
\vol 252 \yr 2002 \pages 27--64
\endref

\ref\key{{\bf 15}}
\by R.M. Green
\paper On planar algebras arising from hypergroups
\jour J. Algebra
\vol 263 \yr 2003 \pages 126--150
\endref

\ref\key{{\bf 16}}
\by R.M. Green
\paper Categories arising from tabular algebras
\jour Glasgow Math. J.
\vol 45 \yr 2003 \pages 333--352
\endref

\ref\key{{\bf 17}}
\by R.M. Green
\paper Standard modules for tabular algebras
\jour Algebr. Represent. Theory
\vol 7 \yr 2004 \pages 419--440
\endref

\ref\key{{\bf 18}}
\by P. Hanlon and D. Wales
\paper A tower construction for the radical in Brauer's centralizer algebras
\jour J. Algebra
\vol 164 \yr 1994 \pages 773-830
\endref

\ref\key{{\bf 19}}
\by V.F.R. Jones
\paper A quotient of the affine Hecke algebra in the Brauer algebra
\jour L'En\-seigne\-ment Math. 
\vol 40 \yr 1994 \pages 313--344
\endref

\ref\key{{\bf 20}}
\by V.F.R. Jones
\paper Planar Algebras, I
\miscnote preprint
\endref

\ref\key{{\bf 21}}
\by D. Kazhdan and G. Lusztig
\paper Representations of Coxeter groups and Hecke algebras
\jour Invent. Math. 
\vol 53 \yr 1979 \pages 165--184
\endref

\ref\key{{\bf 22}}
\by S. K\"onig and C.C. Xi
\paper Cellular algebras: inflations and Morita equivalences
\jour Jour. L.M.S.
\vol 60 \yr 1999 \pages 700--722
\endref

\ref\key{{\bf 23}}
\by M. Kosuda
\paper Representation of $q$-analogue of rational Brauer algebras
\jour Tsukuba J. Math 
\vol 21 \yr 1997 \pages 707--728
\endref

\ref\key{{\bf 24}}
\by G. Lusztig
\paper Cells in affine Weyl groups, II
\jour J. Algebra
\vol 109
\yr 1987
\pages 536--548
\endref

\ref\key{{\bf 25}}
\by G. Lusztig
\paper Cells in affine Weyl groups, IV
\jour J. Fac. Sci. Tokyo U. (IA)
\vol 36
\yr 1989
\pages 297--328
\endref

\ref\key{{\bf 26}}
\by P.P. Martin
\book Potts Models and related problems in statistical mechanics
\publ World Scientific
\publaddr Singapore
\yr 1991
\endref

\ref\key{{\bf 27}}
\by P.P. Martin
\paper On Schur--Weyl duality, $A_n$ Hecke algebras and quantum $sl(N)$
\jour Int. J. Mod. Phys. A
\vol 7 \yr 1992 \pages 645--674
\endref

\ref\key{{\bf 28}}
\by P.P. Martin
\paper Temperley--Lieb algebras for nonplanar statistical
mechanics---the partition algebra construction
\jour J. Knot Th. Ram.
\vol 3 \yr 1994 \pages 51--82
\endref

\ref\key{{\bf 29}}
\by P.P. Martin
\paper The structure of the partition algebras
\jour J. Algebra 
\vol 183 \yr 1996 \pages 319--358
\endref

\ref\key{{\bf 30}}
\by P. Martin and H. Saleur
\paper The blob algebra and the periodic Temperley--Lieb algebra
\jour Lett. Math. Phys.
\vol 30 (3)
\yr 1994 
\pages 189--206
\endref

\ref\key{{\bf 31}}
\by H. Rui and C.C. Xi
\paper The representation theory of cyclotomic Temperley--Lieb algebras
\jour Comment. Math. Helv.
\vol 79 \yr 2004 \pages 427--450
\endref

\ref\key{{\bf 32}}
\by R.P. Stanley
\paper Subdivisions and local $h$-vectors
\jour J. Amer. Math. Soc.
\vol 5 \yr 1992 \pages 805--851
\endref

\ref\key{{\bf 33}}
\by V.S. Sunder
\paper $\text{II}_1$ factors, their bimodules and hypergroups
\jour Trans. Amer. Math. Soc.
\vol 330 \yr 1992 \pages 227--256
\endref

\ref\key{{\bf 34}}
\by H.N.V. Temperley and E.H. Lieb
\paper Relations between percolation
and colouring problems and other graph theoretical problems associated
with regular planar lattices: some exact results for the percolation
problem
\jour Proc. Roy. Soc. London Ser. A 
\vol 322 \yr 1971 \pages 251--280
\endref

\ref\key{{\bf 35}}
\by H. Wenzl
\paper On the structure of Brauer's centralizer algebras
\jour Ann. of Math.
\vol 128
\yr 1988
\pages 173--193
\endref

\ref\key{{\bf 36}}
\by C.C. Xi
\paper Partition algebras are cellular
\jour Compositio Math.
\vol 119 \yr 1999 \pages 99--109
\endref

\ref\key{{\bf 37}}
\by N. Xi
\paper The based ring of two-sided cells of affine Weyl groups of type 
$\tilde A_{n-1}$
\jour Mem. Amer. Math. Soc.
\vol 157 \yr 2002 \pages no. 749
\endref

\endRefs
